\documentclass[leqno,11pt]{article}
\usepackage{latexsym}
\usepackage{amssymb}
\usepackage{amsfonts}
\usepackage{amsmath}
\usepackage{color}
\usepackage[T1]{fontenc}

\newtheorem{definition}{Definition}[section]
\newtheorem{lemma}[definition]{Lemma}
\newtheorem{theorem}[definition]{Theorem}
\newtheorem{proposition}[definition]{Proposition}
\newtheorem{corollary}[definition]{Corollary}
\newtheorem{remark}[definition]{Remark}
\newtheorem{example}[definition]{Example}

\makeatletter
\long\def\unmarkedfootnote#1{{\long\def\@makefntext##1{##1}\footnotetext{#1}}}
\makeatother

\def\w0{{W_0^{1,A}(\Omega)}}
\def\r{\mathbb R}

\def\ep{\varepsilon}

\def\rn{{{\r}^n}}

\def\abar{\mathcal{A}}

\newcommand{\medint}{-\kern  -,395cm\int}
\newcommand{\medintinrigo}{-\kern  -,315cm\int}
\newcommand{\medelle}{-\kern  -,235cm L}
\newcommand{\medellenrigo}{-\kern  -,180cm L}
\newcommand{\qed}{\thinspace\null\nobreak\hfill
\hbox{\vbox{\kern-.2pt\hrule height.2pt
depth.2pt\kern-.2pt\kern-.2pt \hbox to1.8mm {\kern-.2pt\vrule
width.4pt \kern-.2pt\raise1.8mm\vbox to.2pt{} \lower0pt\vtop
to.2pt{}\hfil\kern-.2pt \vrule
width.4pt\kern-.2pt}\kern-.2pt\kern-.2pt \hrule height.2pt
depth.2pt \kern-.2pt}}\par\medbreak}

\title{Existence and regularity results for nonlinear elliptic equations 
in Orlicz spaces}  \numberwithin{equation}{section}

\author{
 Giuseppina Barletta\\
 {\it Dipartimento di Ingegneria Civile, dell'Energia, dell'Ambiente e dei Materiali}
\\ {\it Universit\`a Mediterranea di Reggio Calabria}
 \\ {\it Via Zehender, 89122, Reggio Calabria
 (Italy)}
 \\ {\it e-mail: giuseppina.barletta@unirc.it}
}

\date{}

\begin{document}
	
\maketitle

\begin{abstract} \noindent
We are concerned with the existence and regularity of the solutions to the Dirichlet problem, for a class of quasilinear elliptic equations driven by a general differential operator, depending on $(x,u,\nabla u)$, and with a convective term $f$. The assumptions on the members of the equation are formulated in terms of Young's functions, therefore we work in the Orlicz-Sobolev spaces. After establishing some auxiliary properties, that seem new in our context, we present two existence and two regularity results. We conclude with several examples. 

\end{abstract}

\unmarkedfootnote {
	\smallskip
\par\noindent {\it
The author is member of the Gruppo Nazionale  per l'Analisi Matematica, la Probabilit\`{a} e le loro Applicazioni  (GNAMPA) of the Istituto Nazionale di Alta Matematica (INdAM).\\
 This research is partially supported by the Ministry of Education, University and Research of Italy, Prin 2017 Nonlinear Differential Problems via Variational, Topological and Set-valued Methods (Project No. 2017AYM8XW).
}}

\section{Introduction}\label{sec1}

In the present paper we deal with the existence and regularity of solutions to the
the following nonlinear elliptic problem 
\begin{equation}\label{problem}
	\begin{cases}
		- {\rm div}(\abar (x,u,\nabla u)) =f(x,u, \nabla u) & {\rm in}\,\,\, \Omega \\
		u=0 & {\rm on}\,\,\,\partial \Omega \,.
	\end{cases}
\end{equation}
$\Omega$ is a subset of $\rn$ having finite measure, $\abar :\Omega \times \r\times \rn \to \rn$ 
and $f:\Omega\times\r\times\rn\to\r$ are Carath\'{e}odory functions. \\
In spite of the many boundedness (see, for instance, \cite{BaCiMa, Cbound, LU, MaWi} and the references therein)  or regularity (see \cite {BaCoMi, BPS, DBe, Li, Li1,  Mi}) results, in which the existence of a solution is assumed a priori, there are fewer results about the existence of solutions for \eqref{problem}, unless the principal part is the $p-$Laplacian or growths like the $p-$Laplacian.\\
This is due to the fact that when the operator $\abar$ in \eqref{problem} depends also from the unknown $u$, technical difficulties arise. First of all, variational methods are not directly  applicable (and this holds also for problems with a convective term, namely with a nonlinearity $f$ depending also on $\nabla u$), and the properties of $\abar$ are actually well known only in standard contexts (see \cite{CLM, Tr}). This means, in particular, that in most situations, the abstract framework for the study of \eqref{problem} is the classical Sobolev space $W_0^{1,p}(\Omega)$.\\
In this regard, we note that there are many papers dealing with the existence of solutions to \eqref{problem}, whenever the principal part does not depends explicitely on $u$, and the tecniques used are the most diverse. This also testifies to the growing interest in recent years for the study of these problems.\\
\noindent The two main motivations of the paper are to investigate problems whose differential part depends explicitly also on the unknown, and to remove restrictions of power type growth in the study of \eqref{problem}. To do this, we have to establish some properties of $\abar$, which are of independent interest (see Proposition \ref{S1properties}). Just to give an idea, in our contest we can manage problems driven by an operator having a power-times-logarithmic type growth, like (see Example \ref{exA})
\begin{align}\label{Aintro}
\abar(x,s,\xi)=&
a(x)|s|^\beta\lg^{\beta_1}(1+|s|)|\xi|^{p-2-\delta}\lg^{q(1-\frac\delta{p-1})}(1+|\xi|)\xi+|\xi|^{p-2}\lg^{q}(1+|\xi|)\xi  \,.
\end{align}
To the best of our knowledge no results are actually available for problems driven by such operators.
Conversely, many existence and regularity results are available for \eqref{problem}, when $\abar$ is the $p$-Laplacian, or the $(p,q)$-Laplacian, namely when $\abar(\nabla u)=|\nabla u|^{p-2}\nabla u$ or $\abar(\nabla u)=|\nabla u|^{p-2}\nabla u+ |\nabla u|^{q-2}\nabla u$ (see \cite{AAp, AF, BalFil, FMP, LMZ, MotWin, NgSch, Ruiz, Zou}). There is also an extensive literature concerning problems in which the structure of $\abar$ allows to tackle \eqref{problem} with the same techniques adopted for the $p-$Laplacian.  We refer to \cite{FM, PZ} for problems with $\abar(\nabla u)$, to \cite{Boc, Boc1, BNV, HT, MotWin, NgSch, Tanaka} for problems with $\abar(x,\nabla u)$ and finally, for problems with $\abar(x,u,\nabla u)$ in $W_0^{1,p}(\Omega)$, we cite \cite{FMMT, G}.\\
 Much less is available regarding more general operators, built via a function radial with respect to $\xi$, but not necesseraly a polinomyal. This new situation requires the use of Young functions and Orlicz spaces. Existence (and regularity) results for problems with $\abar(\nabla u)=a(|\nabla u|)\nabla u$ can be found in \cite{BaTo1, BaTo2, CGSS, FMST}. In \cite{DF} the authors deal with an operator depending on the three variables via Young functions of a real variable.\\
 In this paper, the growth conditions on the terms appearing in \eqref{problem} require to replace the customary Sobolev space with an Orlicz space (see \eqref{Aintro}). These conditions cover several instances already studied in the papers cited above. We stress that Young's functions are also involved in the growth of the convective term $f$. Similar hypotheses can be found in \cite{BaTo1, BaTo2, DF}.\\
 Given the non-variational nature of the problem, we use the method of sub and super solutions, togheter with truncation techniques and the theorem of existence of zeros for monotone operators. For the regularity results, a main tool is a Theorem due to Lieberman (see  \cite[Theorem 1.7]{Li1} and Proposition \ref{Lib}).\\
The first step necessary to obtain our results consists in establishing some properties of the operator $\abar$ in Orlicz spaces. This properties, as well as basic definitions and some other auxiliary results, are collected in Section \ref{sec2}. The main existence theorems (and a general example of function $\abar$ satisfying our assumptions) are presented in Section \ref{sec3}: beside to a general existence result, Theorem \ref{Tnew}, in Theorem \ref{Teo1} we consider a special instance of \eqref{problem}, where a subsolution and a supersolution are obtained via variational methods. In Section \ref{sec4} we prove our regularity results. Theorems \ref{Teo1} and \ref{reg} have among the main hypothesis, the existence of a subsolution and a supersolution. Starting from them, we construct a suitable functional, in which appropriate truncations of $\abar $ and of the convective term $f$ are involved.\\
Finally, Section \ref{sec5} is devoted to some examples where we highlight how the method of sub and super solutions works well in all situations where the convective term $f$ has two zeros of opposite sign, namely $f(x,s_1,0)=f(x,s_2,0)=0$ for all $x\in\Omega$ and for some $s_1,s_2\in \r$, enjoying the condition $s_1\cdot s_2<0$. \\
In this regard, it should be noted that the use of the method of sub and super solutions, combined with truncation techniques, is extremely valid when a mix of conditions related to both $\abar$ and $f$ occur. This is due to the fact that once we have a sub and a super solution, we restrict our attention to a suitable truncation of $\abar(x,u,\nabla u)$. If the structure of $f$ allow for sub and super solutions having good properties, such as, for instance, boundedness, then the truncated operator satisfies the hypotheses of the abstract result, even if the growth of $\abar$ is more general than one would expect (see Examples \ref{ex1} and \ref{ex3}).

\section{Preliminaries}\label{sec2}

In this Section we give the main definitions on Young functions and introduce the Orlicz Sobolev spaces that we use in the sequel. We also collect some new results, auxiliary for the proof of the main theorems. Classical results concerning Young functions and Orlicz spaces can be found in \cite{Chl, KrRu, RR1, RR2}. 
\begin{definition} A function $A: [0, \infty ) \to [0, \infty]$ is called a Young function if it is convex, vanishes at $0$, and is neither identically equal to $0$, nor to infinity.
\end{definition}
It is not restrictive to assume that any finite valued Young function is continuous.
For Young functions 
\begin{equation}\label{convexity}
A(\lambda t) \leq \lambda  A(t)\quad  \hbox{for  $\lambda \leq 1$ and $t \geq 0$.}
\end{equation}
\begin{definition}\label{conjugate} The Young conjugate of a Young function $A$ is the Young function $\widetilde A$ defined as
$$\widetilde A (s) = \sup\{ st - A(t):\ t \geq 0\} \quad \hbox{for $s \geq 0$.}$$
\end{definition}
The following inequalities are a consequence of Definition \ref{conjugate}
\begin{equation}\label{A1}
	\frac{A(t)}{t}\leq \widetilde A^{-1}(A(t))\leq 2\frac{A(t)}{t}\quad \hbox{for}\ t> 0\,,
\end{equation}
\begin{equation}\label{A2}
st\leq \widetilde A(s)+A(t)\quad \hbox{for}\ s,t> 0\,.
\end{equation}
The inverse function in \eqref{A1} is the generalized right continuous inverse. In general a Young function is left continuous (in effect it is continuous, unless it takes the value $+\infty$), and when we deal with the inverse we consider the right continuous one.
\begin{definition}
 A Young function $A$ is said to dominate another Young function $B$ near infinity, if there exist constants $k>0$ and $M\geq 0$ such that \begin{equation}\label{dyoung}
   B(t)\leq A(kt)\ \ \hbox{if }\, t\geq M\,.
\end{equation}
If \eqref{dyoung} holds with $M=0$, then we say that $A$ dominates $B$ 
 Two Young functions $A$ and $B$ are called equivalent near infinity (globally) if they dominate each other near infinity (globally) and w
 e briefly write $A\approx B$ near infinity (globally).
\end{definition}
\begin{definition} 
Two  functions $f,g : (0, \infty)\to [0, \infty )$, are equivalent near infinity (briefly $f \approx g $ near infinity) if and only if there exist suitable positive constants $c_1$, $c_2$ and $s_0$ such that
 \begin{equation}\label{approx1}
 	c_1g(c_1s) \leq f(s) \leq c_2g(c_2s)\quad \hbox{if} \ s>s_0\,.
 \end{equation}
\end{definition}
\begin{definition} A Young function $B$ is said to increase essentially more slowly than $A$ near infinity (briefly $B\ll A$), if $B$ is finite valued and
	\begin{equation}\label{<<}
		\lim_{t\to \infty}\frac{B(\lambda t)}{A(t)}=0\ \ \hbox{for all}\; \lambda >0\,.
	\end{equation}
\end{definition}
Condition \eqref{<<} is equivalent to 
\begin{equation}\label{<<-1}
	\lim_{t\to \infty}\frac{A^{-1}(t)}{B^{-1}(t)}=0\,,
\end{equation}
where $A^{-1}$ denotes the right continuous inverse of $A$. Note that if $A(t)=+\infty$ for $t>t_0$, then $A^{-1}(s)$ is constant for $s$ sufficiently large, whatever inverse we take.
\begin{proposition}\label{propdom}
	Let $A$ and $B$ be two Young functions, such that $B\ll A$. Then there exists a Young function $C$ such that $B\ll C\ll A$.
\end{proposition}
{\bf Proof. } Let us define $C^{-1}(t)=\sqrt{A^{-1}(t)B^{-1}(t)}$ for $t\geq 0$. One can easily check that $C^{-1}$ is concave, strictly increasing (inasmuch $A$ and $B$ are strictly increasing), and 
\begin{align}
	\lim_{t\to \infty}\frac{A^{-1}(t)}{C^{-1}(t)}=0,\ &\ \lim_{t\to \infty}\frac{C^{-1}(t)}{B^{-1}(t)}=0\,. 
\end{align}
Thus $C$ satisfies $B\ll C\ll A$.\qed
\begin{definition} A  Young function $A$ is said to satisfy the $\Delta_2$-condition near infinity (briefly $A\in \Delta_2$ near infinity) if it is finite valued and there exist two constants $K\geq 2$ and $M\geq 0$ such that
\begin{equation}\label{delta2young}
A(2t)\leq KA(t)\quad \hbox{for }\ t\geq M\,.
\end{equation}
\end{definition}
\begin{definition} The function $A$ is said to satisfy the $\nabla_2$-condition near infinity (briefly $A\in\nabla_2$ near infinity) if there exist two constants $K>2$ and  $M\geq 0$ such that
\begin{equation}\label{nabla2young}
A(2t)\geq KA(t)\quad \hbox{for }\ t\geq M\,.
\end{equation}
\end{definition}
If \eqref{delta2young} or \eqref{nabla2young} holds with $M=0$, then $A$ is said to satisfy the $\Delta_2$-condition (globally), or the $\nabla_2$-condition (globally), respectively. If \eqref{delta2young} or \eqref{nabla2young} holds for $0\leq t\leq M$, then $A$ is said to satisfy the $\Delta_2$-condition near zero, or the $\nabla_2$-condition near zero, respectively. 
If $A\in\Delta_2\cap \nabla_2$ near infinity then $\widetilde{A}\in\Delta_2\cap \nabla_2$ near infinity (see \cite{BaCi1}, Proposition 6.6).\\

\par
We give basic definitions and the main properties on the Orlicz spaces.
Let $\Omega$ be a measurable set in $\rn$, with $n\geq 1$. Given a Young function $A$, the Orlicz space $L^A(\Omega)$ is the set of all measurable functions $u:\Omega\to\r$ such that the Luxemburg norm
$$\|u\|_{L^A(\Omega)}=\inf\bigg\{\lambda>0:\int_\Omega A\big(\tfrac 1\lambda|u|\big)\,dx\leq 1 \bigg\}$$
is finite. The functional $\|\cdot\|_{L^A(\Omega)}$ is a norm on $L^A(\Omega)$, and it is a Banach space (see \cite{Adams}). If $A$ and $B$ are Young functions, $A\approx B$ near infinity and $|\Omega|<\infty$ then $L^A(\Omega)=L^B(\Omega)$, and there exist $k_1(,t_0,|\Omega|), k_2(t_0,|\Omega|)$ such that
\begin{align}\label{normequiv}
	k_1\|u\|_{L^A(\Omega)}\leq \|u\|_{L^B(\Omega)}\leq k_2\|u\|_{L^A(\Omega)},\ \hbox{for all}\ u\in L^A(\Omega)\,,\end{align}
If $A$ is a Young function and $\widetilde A$ denotes its conjugate, then a generalized  H\"older inequality
\begin{equation}\label{holderyoung}
\int _\Omega |u v|\,dx \leq 2\|u\|_{L^A (\Omega)} \|v\|_{L^{\widetilde A}(\Omega)}
\end{equation}
holds for every $u\in L^A (\Omega)$ and $v\in L^{\widetilde A}(\Omega)$.\\
If $A\in \Delta_2$ globally (or $A\in \Delta_2$ near infinity and $\Omega$ has finite measure) then 
\begin{equation}\label{intfinito}
 \int_{\Omega} A(k|u|)dx <+\infty \ \hbox{for all}\ u\in L^A(\Omega),\ \hbox{all}\ k\geq 0\,. 
\end{equation}
Also, if $B\ll A$ near infinity and $\Omega$ has finite measure, then 
\begin{equation}\label{intfinito1}
	\int_{\Omega} B(k|u|)dx <+\infty \ \hbox{for all}\ u\in L^A(\Omega),\ \hbox{all}\ k\geq 0\,. 
\end{equation}
Given a Young function $A \in C^1([0,+\infty))$, define the quantities 
\begin{align}\label{i}
i_A=\inf_{t>0}\frac {t\cdot A'(t)}{A(t)},\,s_A=\sup_{t>0}\frac {t\cdot A'(t)}{A(t)},\ &\ i_A^\infty=\liminf_{t\to +\infty}\frac {t\cdot A'(t)}{A(t)},\, s_A^\infty=\limsup_{t\to +\infty}\frac {t\cdot A'(t)}{A(t)}\,.
\end{align}
The conditions
\begin{align*}\label{ndg}
	i_A>1\ \hbox{and}\ s_A<+\infty\ &\ (i_A^\infty>1\ \hbox{and}\ s_A^\infty<+\infty )
\end{align*}
are equivalent to the fact that $A\in\nabla_2\cap\Delta_2$ globally (near infinity). The following result extends Lemma A.1 of \cite{CM} to Young functions $A\in\nabla_2$ near infinity. 
\begin{lemma}
Let $\Omega$ be a subset of $\rn$, with finite measure. Let $A$ be a Young function. Assume that $i_A^\infty>1$. Then there exists $k_3(i_A^\infty,\,|\Omega|)\geq 0$
 such that 
\begin{equation}\label{min per coerc}
	\int_\Omega A(|v|)dx\geq k_1^{i^\infty_A}\|v\|_{L^A(\Omega)}^{i^\infty_A}-k_3 \quad  \hbox{for every}\ v\in L^A(\Omega), \ \hbox{such that}\ \|v\|_{L^A(\Omega)}\geq \frac 1{k_1}\,.
\end{equation}
Here $k_1$ is that of \eqref{normequiv}.
\end{lemma}	
{\bf Proof}. Since $i_A^\infty>1$, corresponding to $\ep\in (0,i_A^\infty-1) $ there exists $t_0>0$ such that $\displaystyle{\frac {t\cdot A'(t)}{A(t)}>i_A^\infty-\ep}$ for $t\geq t_0$. 
Fix $\alpha>i_A^\infty$ and consider the Young function
\begin{align}
	A_1(t)=\left\{\begin{array}{ll} c_1t^{i_A^\infty}+c_2t^\alpha & \hbox{if}\ t\leq t_0\\
		A(t)&\hbox{if}\ t>t_0\,	
\end{array}\right.
\end{align}
where $c_1$ and $c_2$ are chosen in such a way that $A_1\in C^1([0,+\infty[)$. Then $A_1\approx A$ near infinity and $i_{A_1}^\infty=i_{A}^\infty$. Let $v\in L^A(\Omega)$ be such that $\|v\|_{L^A(\Omega)}\geq \frac1{k_1}. $

Using Lemma A.1 of \cite{CM} and \eqref{normequiv}, with $B=A_1$
\begin{align*}
\int_\Omega A(|v|)dx&=\int_\Omega A_1(|v|)dx+\int_{\{|v|\leq t_0\}} 
(A(|v|)-A_1(|v|))dx\\
&\geq \|v\|_{L^{A_1}(\Omega)}^{i_A^\infty}+m_0|\Omega|\geq
k_1^{i_A^\infty}\|v\|_{L^{A}(\Omega)}^{i_A^\infty}+ 
m_0|\Omega|\,,
\end{align*}
where $m_0=\min_{0\leq t\leq t_0}(A(t)-A_1(t))\leq 0$. Thus \eqref{min per coerc} follows with 
$k_3= -m_0|\Omega|\geq 0$. \qed
The Lemma below is a consequence of the Theorems of Vitali and De la Vall\'ee- Poussin. It will be used several times in the paper.
\begin{lemma}\label{corvit}
Let $\Omega\subseteq\rn$ be such that $|\Omega|<\infty$ and let $A$ be a Young function. 
If $\{u_k\}_{k\in N}\subseteq L^A(\Omega)$ satisfies
\begin{itemize}
\item $u_k(x)\to u(x) $ a.e. in $\Omega$\,, 
\item $\{u_k\}_{n\in N}$ is bounded in $L^A(\Omega)$\,,
	\end{itemize}
then $u_k\to u$ in $L^B(\Omega)$ for all Young functions $B\ll A$ near infinity.\\
If $L^A(\Omega)$ is reflexive, then $|u_k-u|\rightharpoonup 0$ in $L^A(\Omega)$.
\end{lemma}	
{\bf Proof.} Let $B\ll A$ near infinity. Then $B$ is finite valued and we may assume, without loss of generality, that $B$ is continuous. We must prove that $\lim_{k\to +\infty}\int_\Omega B\left(\frac{|u_k-u|}{\lambda}\right)dx=0$, for all $\lambda >0$.\\
If $A(t)\equiv +\infty$ for $t>t_0>0$, then $L^A(\Omega)=L^\infty(\Omega)$ and the conclusion follows via the dominated convergence theorem.\\
If $A$ is finite valued, then we may assume that it is continuous. Fix $\lambda >0$ and let $M>0$ be such that $\|u_k\|_{L^A(\Omega)}\leq M$ for all $n\in N$. Then, the continuity of $A$ and the Fatou's lemma guarantee that $u\in L^A(\Omega)$ and $\|u\|_{L^A(\Omega)}\leq M$. Let $\Psi(t)= A\big(\frac{\lambda}{2M}B^{-1}(t)\big)$, for $t\geq 0$. Then $\Psi$ is increasing and the condition $B\ll A$ guarantees that
\begin{align}
	\lim_{t\to +\infty} \frac{\Psi(t)}{t}=+\infty\,.
\end{align}
Note that 
\begin{align}
	\left\{B\left(\frac{|u_k-u|}{\lambda}\right)\right\}_{k\in N}\subseteq \left\{v\in L^1(\Omega):\,\int_\Omega \Psi(|v|)dx\leq 1\right\}\,.
\end{align} 
Thus the family $\left\{B\left(\frac{|u_k-u|}{\lambda}\right)\right\}_{k\in N}$ is equintegrable and we can apply the Vitali's Theorem. This proves that $u_k\to u$ in $L^B(\Omega)$.\\
For the last part of the proof we observe that $|u_k-u|\rightharpoonup v$
in $L^A(\Omega)$, up to a subsequence. From the proof above we know that $|u_k-u|\to 0$ in $L^B(\Omega)$. Thus $v\equiv 0$.
 This apply to all the subsequences, thus it holds for the whole sequence. \qed

From now and throughout the paper we assume that $\Omega$ is an open set in $\rn$ with $|\Omega|<\infty$. The isotropic Orlicz-Sobolev spaces $W^{1,A}(\Omega)$ and $W^{1,A}_0(\Omega)$ are defined as
\begin{align}\label{W1A}
	W^{1,A}(\Omega)=\{u:\Omega\to\r:&\, u\hbox{ is weakly differentiable in $\Omega$,\ $u$,\;$|\nabla u| \in L^A (\Omega)$}\}\,
\end{align}
and
\begin{align}\label{W01A}
	\w0=\{u:\Omega \to \r: &\, \hbox{the continuation of $u$ by $0$ outside $\Omega$} \\ & \hbox{is weakly differentiable in $\rn$, \; u,\; $|\nabla u| \in L^A (\Omega)$}\}.\nonumber
\end{align}
The spaces  $W^{1,A}(\Omega)$ and $\w0$ equipped with the norms
\begin{align}
	\|u\|_{W^{1,A}(\Omega)} = \|u\|_{L^A(\Omega)}+\|\nabla u\|_{L^A (\Omega)},\ \  \hbox{and}\ \ \|u\|_{\w0} = \|\nabla u\|_{L^A (\Omega)},
\end{align}
are Banach spaces. 
The norm on $\w0$ is equivalent to the standard one
$$\|u\|_{W_0^{1,A}(\Omega)} = \|u\|_{L^A(\Omega)}+\|\nabla u\|_{L^A (\Omega)}\,.$$
If $A\in \Delta_2\cap \nabla_2$ near infinity, then the Orlicz-Sobolev space $W^{1,A}_0(\Omega)$ is reflexive (\cite [Proposition 3.1] {BaCi1}).
\begin{definition}
The optimal Sobolev conjugate of $A$ is defined by $A_n:[0,\infty)\to [0,\infty]$
\begin{equation}\label{sobcong}
A _n (t)= A({H}^{-1}(t)) \quad \hbox{for $t \geq 0$,}
\end{equation}
where ${ H}:[0,\infty ) \to [0, \infty)$ is given by
\begin{equation*}\label{H1}
{ H}(t)=\bigg(\int _0^t\bigg(\frac\tau{A(\tau)}\bigg)^{\frac 1{n-1}}\,d\tau\bigg)^{\frac{n-1}{n}}\quad \hbox{for $t \geq 0$,} \end{equation*}
provided that the integral is convergent.
Here, ${ H}^{-1}$ denotes the generalized left-continuous inverse of ${H}$.
\end{definition}
If 
\begin{equation}\label{conv0}
\int_0\bigg(\frac\tau{A(\tau)}\bigg)^{\frac 1{n-1}}\,d\tau< \infty\,,
\end{equation}
then (see \cite[Theorem 1]{Cfully})
\begin{equation}\label{W0An}
	\w0\to L^{A_n}(\Omega)\,.
\end{equation}
and there exists a constant $C=C(n)$ such that
\begin{equation}\label{B-Wbis}	\|u\|_{L^{A _n}(\Omega )}\leq C\|u\|_{W_0^{1,A}(\Omega)}\quad \hbox{for every}\ u \in W_0^{1,A} (\Omega)\,.\end{equation}
Thus from  \cite[Theorem 3]{Csharp} (or \cite[Theorem 3.1]{HLZ}), the embedding $\w0\to L^{E}(\Omega)$ is compact for every Young function $E\ll A_n$ near infinity.\\
If
\begin{equation}\label{div}
	\int^\infty  \bigg(\frac \tau{A(\tau)}\bigg)^{\frac1{n-1}}\, d\tau = \infty
\end{equation}
then $A_n$ assumes only finite values, while when
\begin{equation}\label{conv}
	\int ^\infty  \bigg(\frac \tau{A(\tau)}\bigg)^{\frac1{n-1}}\, d\tau < \infty\,,
\end{equation}
then $A_n(t)=\infty$ for $t$ large and \eqref{B-Wbis} yields
\begin{equation*}
	\|u \|_{L^{\infty}(\Omega )} \leq C \|u \|_{W_0^{1,A}(\Omega)}\quad \hbox{for every $u\in \w0$.}
\end{equation*}
When \eqref{div} holds then (see \cite[(3.13)]{BaCi1})
\begin{equation}\label{Anfin}
	\int_\Omega A_n(\lambda|u|)dx<\infty\ \hbox{for every}\ u\in W_0^{1,A}(\Omega),\ \hbox{and every}\, \lambda>0.
\end{equation}

\begin{example}{\bf A general Young function}\label{ex1}\\
Let $1<p<+\infty$, $q\in \r$. Consider a Young function $B:[0,+\infty[\rightarrow [0,+\infty[$, such that 
\begin{align}\label{B'ex}
	B'(t)= t^{p-1}\lg^q(1+t)\quad \hbox{for}\ t>>1\,.
\end{align} 
Then 
\begin{align}\label{reflexivity}
	i_{B}^\infty=s_{B}^\infty=p>1\,,
\end{align}
and $W_0^{1,B}(\Omega)$ is reflexive. Also
\begin{align}\label{Bequiv}
B(t)\approx t^{p}\lg^{q}(1+t)\quad \hbox{near infinity}.
\end{align}
Thus, given any Young function $A$ such that 
\begin{align}\label{Aex}
	A(t)\approx t^p\lg^q(1+t)\ \hbox{near infinity}\,,
\end{align} 
and a set $\Omega$ with finite measure, it holds $\w0=W_0^{1,B}(\Omega)$, up to equivalent norms. The conjugate $\widetilde{A}$ and the optimal Sobolev conjugate of $A$ satisfy 
\begin{align}\label{tilde}
	\widetilde{A}(t)\approx t^{p'}\lg^{-\frac q{p-1}}(1+t)\quad \hbox{near infinity},
\end{align}
and
\begin{align}\label{optimal}
	A_n(t) \thickapprox \left\{\begin{array}{lll}
		t^{\frac {np}{n-p}}\lg^{\frac {nq}{n-p}}(1+t)&\hbox{near infinity}&\hbox{when}\ 1<p<n,\\
		e^{t^{\frac {n}{n-q-1}}}&\hbox{near infinity}&\hbox{when}\ p=n,\ q<n-1,\\
		e^{e^{t^{\frac {n}{n-1}}}}&\hbox{near infinity}&\hbox{when}\ p=n,\ q=n-1,\\	
		+\infty&\hbox{near infinity}&\hbox{in all the other cases}\,.
	\end{array}\right.
\end{align}
If $p+q-1\geq 0$ then the function $B$ satisfying \eqref{B'ex} for all $t>0$ is in fact a Young function. When $p+q-1>0$ the function $B\in\nabla _2$ globally.
All the examples of the paper will involve functions $A$ complying with \eqref{Aex}. 
\end{example}

We recall now some definitions and the theorem on pseudomonotone operators. Then we introduce the conditions on the function $\abar$ in \eqref{problem}.
\begin{definition}
	Let $X$ be a real reflexive Banach space. A mapping  ${\cal B}: X\to X^*$ is called
	\begin{itemize}
		\item [(i)] coercive if $\lim_{\|u\|\to\infty}\frac{\langle {\cal B}u,u\rangle}{\|u\|}=+\infty$;
		\item [(ii)] bounded if it maps bounded sets into bounded sets;
		\item [(iii)] pseudomonotone if $u_n\rightharpoonup  u$ and  $\limsup_{k\to +\infty}\langle {\cal B}u_k,u_k-u\rangle\leq 0$ imply that  ${\cal B}u_k\rightharpoonup  {\cal B}u$ and $\langle {\cal B}u_k,u_k\rangle \to \langle {\cal B}u,u\rangle$.
		\item [(iv)] satisfies the $(S)_+$ property if $u_k\rightharpoonup  u$ and  $\limsup_{k\to +\infty}\langle {\cal B}u_k,u_k-u\rangle\leq 0$ imply that $u_k\to u$ (strongly) in $X$.
	\end{itemize}
\end{definition}
\begin{theorem} \label{theorem on pseudomonotone operator}
	{\rm  (see \cite[Theorem 2.99]{CLM})}	
	Let $X$ be a real reflexive Banach space and let ${\cal B}:X\to X^*$ be a bounded, coercive and pseudomonotone operator. Then, for every $b\in X^*$ the equation ${\cal B}x=b$ has at least solution $x\in X$.
\end{theorem}

Let us now introduce the hypotheses on the function $ \abar $ operator and verify the properties of the integral operator associated with it, in a standard way. We note that this is a Leray Lions type operator. The properties of the auxiliary truncated operator will easily follow from those of the non-truncated operator (see Proposition \ref{S1properties}). We extend the results in \cite{CLM, Tr} in several directions: Orlicz spaces are considered, and even in the case of Lebesgue spaces, the hypotheses are slightly more general.\\
Let $\Omega\subset \rn$ be a set with finite measure and let $A$ be a Young function, $A\in \Delta_2\cap\nabla_2$ near infinity. Consider the vector valued function $\abar:\Omega\times \r\times \rn\to\rn$, $\abar =(a_1,\ldots a_n)$, enjoying with the properties that each $a_i(x,s,\xi)$ is a Carath\'{e}odory function, and
\begin{equation}\label{a1}
	|\abar(x,s,\xi)| \leq q(x)+b\left[\widetilde A^{-1}(F(b|s|))+ \widetilde A^{-1}(A(b|\xi|))\right]\ \hbox{for a.e.} \ x\in\Omega,\ \hbox{all}\ (s,\xi)\in\r\times\rn\,.
\end{equation} 
Here $q\in L^{\widetilde A} (\Omega)$, $b>0$, and $F$ is a Young 
function, $F\ll A_n$ near infinity.
\begin{align}\label{a2}
\sum_{i=1}^n\left(a_i(x,s,\xi)-(a_i(x,s,\xi ')\right)\cdot (\xi_i-\xi_i')>0\ &\ \hbox{for a.e.}\, x\in \Omega,\, \hbox{all}\, s\in\r,\, \hbox{all}\,\, \xi,\xi'\in\rn,\, \xi\neq\xi' \,.
\end{align}
\begin{align}\label{a3bis}
	\sum_{i=1}^na_i(x,s,\xi)\cdot \xi_i\geq cA(c|\xi|)-dG(d|s|)-r(x)&\qquad \hbox{for a. e.}\ x\in \Omega,\ \hbox{all}\ (s,\xi)\in\r\times\rn \,.
	\end{align}
Here  $c,d>0$, $G$ is a Young function, $G\ll A_n$ near infinity, and $r\in L^{1}(\Omega)$,.
\begin{remark}
Using the condition $A\in \Delta_2$ near infinity, and $|\Omega|<\infty$, it can be shown that \eqref{a3bis} is equivalent to \begin{equation}\label{a3}
		\sum_{i=1}^na_i(x,s,\xi)\cdot \xi_i\geq cA(|\xi|)-dG(d|s|)-r(x)\qquad \hbox{for a. e.}\ x\in \Omega,\ \hbox{all}\ (s,\xi)\in\r\times\rn \,.
	\end{equation}
We will use \eqref{a3} in the sequel.
\end{remark}

\begin{example}\label{exA}
We present some examples of functions $\abar(x,s,\xi)$ satisfying \eqref{a1}, \eqref{a2} and \eqref{a3bis}.\\
Let $p,q\in \r$, satisfying 
\begin{align}\label{pqageneral}
	1<p\leq n,\ &\ q\leq n-1,\ \ p-1+q>0\,.
\end{align} 
Consider also a Young function $D$, $D(t)\approx t^r\lg^{-\frac q{p-1}}(1+t)$ near infinity, for some $r>1$. \\
When $p<n$ we take
\begin{align}\label{p<n}
	\abar(x,s,\xi)=&
	a(x)|s|^\beta\lg^{\beta_1}(1+|s|)|\xi|^{p-2-\delta}\lg^{q(1-\frac\delta{p-1})}(1+|\xi|)\xi+|\xi|^{p-2}\lg^{q}(1+|\xi|)\xi  \,,\\
	&\hbox{where}\ 0<\delta<p-1<n-1,\ 0<\beta<\frac{n\delta}{n-p},\ \beta+\beta_1>0,\nonumber\\
	&a(x)>0,\ a\in L^D(\Omega),\ \hbox{and}\ r>\frac{np}{n\delta-\beta(n-p)}.\nonumber
\end{align}
When $p=n$, $q<n-1$, we define
\begin{align}\label{p=n} 
	\abar(x,s,\xi)=&a(x)e^{|s|^\beta}|\xi|^{n-2-\delta}\lg^{q(1-\frac\delta{n-1})}(1+|\xi|)\xi+|\xi|^{n-2}\lg^{q}(1+|\xi|)\xi  \,,\\
	&\hbox{where}\ 0<\delta<n-1,\  0<\beta<\frac{n}{n-q-1},\nonumber \\
	&a(x)>0,\ a\in L^D(\Omega),\
	\hbox{and}\ r>\frac{n}{\delta}.\nonumber
\end{align}
Finally, when $p=n$, $q=n-1$ we consider
\begin{align}\label{q=n-1}
	\abar(x,s,\xi)=&a(x)e^{e^{|s|^\beta}}|\xi|^{n-2-\delta}\lg^{n-1-\delta}(1+|\xi|)\xi+|\xi|^{n-2}\lg^{n-1}(1+|\xi|)\xi  \,,\\
	&\hbox{where}\	0<\delta<n-1,\  0<\beta<\frac{n}{n-1},\nonumber\\
	&a(x)>0,\ a\in L^D(\Omega),\  \hbox{and}\ r>\frac{n}{\delta}.\nonumber
\end{align}   

It holds 
\begin{align}\label{modulo1}
	|\abar(x,s,\xi)|\leq a(x)|s|^\beta\lg^{\beta_1}(1+|s|)|\xi|^{p-1-\delta}\lg^{q(1-\frac\delta{p-1})}(1+|\xi|)+|\xi|^{p-1}\lg^{q}(1+|\xi|)  \,,
\end{align}
\begin{align}\label{modulo2}
	|\abar(x,s,\xi)|\leq a(x)e^{|s|^\beta}|\xi|^{n-1-\delta}\lg^{q(1-\frac\delta{n-1})}(1+|\xi|)+|\xi|^{n-1}\lg^{q}(1+|\xi|)  \,,
\end{align}
and
\begin{align}\label{modulo}
	|\abar(x,s,\xi)|\leq a(x)e^{e^{|s|^\beta}}|\xi|^{n-1-\delta}\lg^{n-1-\delta}(1+|\xi|)+|\xi|^{n-1}\lg^{n-1}(1+|\xi|)  \,,
\end{align}
for all $(x,s,\xi)\in \Omega\times\r\times\rn$, respectively.\\
Now we use twice \eqref{A2} in \eqref{modulo1}, \eqref{modulo2} and \eqref{modulo}: first with $B(t)=t^{\frac{p-1}\delta}$ and then with $B(t)=t^{\frac{r\delta}p}$. 
\begin{align}\label{AY1}
	|\abar(x,s,\xi)|&\leq a(x)^{\frac{p-1}{\delta}}|s|^{\beta\frac{p-1}{\delta}}\lg^{\beta_1\frac{p-1}{\delta}}(1+|s|)+2|\xi|^{p-1}\lg^{q}(1+|\xi|) \\
	&\leq a(x)^{r\frac{p-1}{p}}+|s|^{\beta r\frac{p-1}{r\delta-p}}\lg^{\beta_1r\frac{p-1}{r\delta-p}}(1+|s|)+2|\xi|^{p-1}\lg^{q}(1+|\xi|)\,. \nonumber
\end{align}
\begin{align}\label{AY2}
	|\abar(x,s,\xi)|&\leq a(x)^{\frac{n-1}{\delta}} e^{\frac{n-1}{\delta}|s|^{\beta}}+2|\xi|^{n-1}\lg^{q}(1+|\xi|) \\
	&\leq a(x)^{r\frac{n-1}{n}}+ e^{\frac{r(n-1)}{r\delta-n}|s|^{\beta}}+2|\xi|^{n-1}\lg^{q}(1+|\xi|)\,. \nonumber
\end{align}
\begin{align}\label{AY3}
	|\abar(x,s,\xi)|&\leq a(x)^{\frac{n-1}{\delta}} e^{\frac{n-1}{\delta}e^{|s|^{\beta}}}+2|\xi|^{n-1}\lg^{n-1}(1+|\xi|) \\
	&\leq a(x)^{r\frac{n-1}{n}}+ e^{\frac{(r)n-1}{r\delta-n\delta}e^{|s|^{\beta}}}+2|\xi|^{n-1}\lg^{n-1}(1+|\xi|)\,. \nonumber
\end{align}
Let $A(t)$ be defined as in \eqref{Aex}.
Taking into account the conditions in \eqref{p<n} and in \eqref{p=n}, for the function $a$ in \eqref{AY1} and \eqref{AY2} it holds $a(x)^{r\frac{p-1}{\delta}}\in L^{\widetilde{A}}(\Omega)$. When \eqref{q=n-1} is in force, then the function $a$ in \eqref{AY3} satisfies $a(x)^{r\frac{n-1}{n}}\in L^{\widetilde{A}}(\Omega)$. Also, $\widetilde{A}^{-1}(A(|\xi|))\approx |\xi|^{p-1}\lg^{q}(1+|\xi|)$ near infinity.\\ 
When \eqref{p<n} is in force we consider the Young function 
$F(s)\approx |s|^{\frac{\beta rp}{r\delta-p}}\lg^{\frac{\beta_1rp}{r\delta-p}-\frac q{p-1}}(1+|s|)$ near infinity. 
When \eqref{p=n} holds we choose the Young function 
$F(s)\approx e^{s^{\beta}}$ near infinity. 
Similarly, when \eqref{q=n-1} is in effect, then we choose $F(s)\approx e^{e^{s^{\beta}}}$ near infinity. In all the three cases $F\ll A_n$ near infinity and $\widetilde{A}^{-1}(F(|s|))$ is equivalent, near infinity, to the function in the $s$ variable, appearing in \eqref{AY1}, \eqref{AY2} and \eqref{AY3} respectively.
Then \eqref{a1} holds with $q(x)=a(x)+k$, for some $k>0$.\\
Condition \eqref{a2} is satisfied, because \eqref{pqageneral} guarantees that the functions $B(t)=t^{p-1-\delta}\lg^{q(1-\frac\delta{p-1})}(1+t)$ and $B(t)=t^{n-1-\delta}\lg^{n-1-\delta}(1+t)$ are increasing, for $t\geq 0$. Condition \eqref{a3bis} clearly holds true.\\
We now consider the borderline instances concernig $\beta$. When $\beta=0$, then we can take $r=\frac{p}\delta$. When $\beta=\frac{n\delta}{n-p}$ in \eqref{p<n}, then in addition to $\beta_1>-\frac{n\delta}{n-p}$, it is necessary $\beta_1<\frac{q\delta(n-1)}{(n-p)(p-1)}$ and $a\in L^{\infty}(\Omega)$. \\
When $q(x)\equiv 0$ then we have the classical $p$-Laplacian. 

\end{example}

Define now  the operator ${\cal S}: \w0\to(\w0)^*$ as
\begin{equation}\label{Sdef}
	\langle{\cal S}u,v\rangle=\int _\Omega \abar(x,u,\nabla u)\cdot\nabla v \,dx\quad
	 \hbox{for $u, v \in \w0$.}
\end{equation}

The properties of ${\cal S}$ are listed in the next proposition.
\begin{proposition}\label{S1properties} Let $A$
	be a Young function, $A\in\nabla_2\cap\Delta_2$ near infinity. Assume that $\abar$ satisfies \eqref{a1}, \eqref{a2}, \eqref{a3bis}. Then the operator ${\cal S}$, introduced in \eqref{Sdef}, is well defined, bounded, continuous, and satisfies the $(S_+)$ property.
\end{proposition}
{\bf Proof}. Let $u\in \w0$. From \eqref{a1}, \eqref{Anfin} and Propositon 6.6 of \cite{BaCi1} we can find $c_1>0$ such that
\begin{align}\label{well}
\int_\Omega \widetilde A(|\abar(x,u,\nabla u)|)dx&\leq \frac{1}{3}\int_\Omega \widetilde A(3q(x))dx+\frac{1}{3}\int_\Omega\widetilde A(3b\widetilde A^{-1}(F(b| u|)
))dx +\frac{1}{3}\int_\Omega\widetilde A(3b\widetilde A^{-1}(A(b|\nabla u|))) dx\nonumber\\
	&\leq c_1\left[ \int_\Omega \left( \widetilde A
	(3q(x))+F(b|u|)+A(b|\nabla u|)\right)dx+1\right]=M<+\infty \,. 
\end{align}
Thus $|\abar(x,u,\nabla u)|\in L^{\widetilde A}(\Omega)$ and $\|\abar(x,u,\nabla u) \|_{L^{\widetilde A}(\Omega)}\leq \max\{1,\,M\}$. From \eqref{holderyoung}
\begin{align*}
	|\langle {\cal S}u,v\rangle|\leq 2\|\abar(x,u,\nabla u) \|_{L^{\widetilde A} (\Omega)}\|\nabla v\|_{L^A(\Omega)}\,.
\end{align*}

Similarly, if ${\cal C}\subseteq W^{1,A}_0(\Omega)$ is such that $\|u\|_{W^{1,A}_0(\Omega)}\leq M_0$ for some $M_0>0$, then $\|u\|_{L^{A_n}(\Omega)}\leq M_1$ for all $u\in {\cal C}$. From \eqref{well} and Lemma 2.7 of \cite{BaTo1}, there exist $c_1, L>0$ such that 
\begin{equation*}
	\int_\Omega\widetilde A\left( \abar(x,u,\nabla u)\right)dx\leq
	c_1\left(\int_\Omega \big(\widetilde A
	\left(3q(x)\right)+ F \left(b|u|\right)+A \left(b|\nabla u|\right)\big)\,dx+1\right)\leq L\ \hbox{for all}\ u\in {\cal C}.
\end{equation*}
Then
\begin{equation}\label{abound}
	\left\|\abar(x,u,\nabla u)\right\|_{L^{\widetilde A}(\Omega)}\leq \max\{1,L\}\ \hbox{for all}\ u\in {\cal C}\,,
\end{equation}
and
\begin{equation*}
	\left\|{\cal S}(u)\right\|_{\left(\w0\right)^*}\leq 2\max\{1,L\}\ \hbox{for all}\ u\in {\cal C}\,.
\end{equation*}
Let now $\{u_k\}$ be a sequence in $\w0$ converging to $u\in \w0$. Then, we can find a subsequence, say still $\{u_k\}$, and two functions $g_1\in L^{A_n}(\Omega)$, $g_2\in L^A(\Omega)$ such that $|u_k(x)|\leq g_1(x)$, $|u(x)|\leq g_1(x)$, $|\nabla u_k(x)|\leq g_2(x)$ and $|\nabla u(x)|\leq g_2(x)$ for a.a. $x\in \Omega$ and for all $k\in N$. Let $\lambda >0$. Using \eqref{a1} and the condition $\widetilde A\in \Delta_2$ near infinity we can find $c_1>0$ such that
\begin{align}
\widetilde A\left( \frac{|\abar(x,u_k,\nabla u_k)-\abar(x,u,\nabla u)|}{\lambda}\right)\leq&\frac{1}{5}\widetilde A
\left(\frac{10q(x)}{\lambda}\right)\\
+\frac{1}{5}\bigg[\widetilde A\left(\frac{5b}{\lambda}\widetilde A^{-1}(F \left(b|u_k|\right))\right) 
+\widetilde A\left(\frac{5b}{\lambda}\widetilde A^{-1}(F \left(b|u|\right))\right)+&\widetilde A\left(\frac{5b}{\lambda}\widetilde A^{-1}(A\left(b|\nabla u_k|\right))\right)+\widetilde A\left(\frac{5b}{\lambda}\widetilde A^{-1}( A \left(b|\nabla u|\right))\right)\bigg]\nonumber\\
\leq\frac{1}{5}\widetilde A
\left(\frac{10q(x)}{\lambda}\right)+c_1(1+F \left(b|u_k|\right)	+F \left(b|u|\right)+&A\left(b|\nabla u_k|\right)+A\left(b|\nabla u|\right))\nonumber\\
\leq\frac{1}{5}\widetilde A
\left(\frac{10q(x)}{\lambda}\right)+c_1(1+2F \left(bg_1(x)\right)	+2A\left(bg_2(x)\right))&:=v_\lambda (x)\,.\nonumber
\end{align}
The function $v_\lambda \in L^1(\Omega)$, inasmuch both $A,\widetilde A\in \Delta_2 $ near inifnity and $F\ll A_n$ near infinity. The continuity of $\abar(x,\cdot, \cdot)$ guarantees that $\lim_{k\to +\infty}\abar (x,u_k(x),\nabla u_k(x))=\abar (x,u(x),\nabla u(x))$ a.e. in $\Omega$. We can thus apply the Lebesgue theorem to obtain 
\begin{align}\label{continuity}
	\lim_{k\to+\infty}\int_\Omega \widetilde A\left( \frac{|\abar(x,u_k,\nabla u_k)-\abar(x,u,\nabla u)|}{\lambda}\right)dx=0\ \hbox{for all}\ \lambda >0\,.
\end{align} 
The arguments above apply to any subsequence of $\{u_k\}$. This means that given any  subsequence, we can extract a subsubsequence for which \eqref{continuity} holds. Thus \eqref{continuity} holds for the whole sequence. Note that 
\begin{align*}
	\left\|{\cal S}(u_k)-{\cal S}(u)\right\|_{\left(\w0\right)^*}=\sup_{\|v\|_{\w0}\leq 1}|\langle{\cal S}(u_k)-{\cal S}(u),v\rangle |\leq 2\|\abar(x,u_k,\nabla u_k)-\abar(x,u,\nabla u)|\|_{L^{\widetilde A}(\Omega)}\,.
\end{align*}
Thus 
\begin{align}
	\lim_{k\to+\infty}\left\|{\cal S}(u_k)-{\cal S}(u)\right\|_{\left(\w0\right)^*}=0\,.
\end{align}
This proves the continuity af ${\cal S}$. \\
Let us now demonstrate the $(S)_+$ property. Let $\{u_k\}$ be  a sequence in $\w0$, $u_k\rightharpoonup  u$ and  
\begin{equation}\label{1}
	\limsup_{k\to +\infty}\langle {\cal S}u_k,u_k-u\rangle\leq 0\,.\end{equation}
We divide the proof in four steps. \\
Step 1 $\nabla u_k\to \nabla u$ a.e.\\
Let $F$ be as in \eqref{a1}. From Proposition \ref{propdom} there exists a Young function $F_1$ such that
$F\ll F_1\ll A_n$ near infinity. Thus $\{u_k\}$ strongly converges to $u$ in $L^{F_1}(\Omega)$, and we can find a function $g\in L^{F_1}(\Omega)$ and a subsequence of $\{u_k\}$, say still $\{u_k\}$, such that $|u_k(x)|,|u(x)|\leq g(x)$ a.e. in $\Omega$. Thus 
\begin{align}
	\widetilde A\left( \frac{|\abar(x,u_k,\nabla u)-\abar(x,u,\nabla u)|}{\lambda}\right)\leq&\frac{1}{5}\widetilde A
	\left(\frac{10q(x)}{\lambda}\right)\\
	+\frac{1}{5}\bigg[\widetilde A\left(\frac{5b}{\lambda}\widetilde A^{-1}(F \left(b|u_k|\right))\right) 
	+\widetilde A\left(\frac{5b}{\lambda}\widetilde A^{-1}(F \left(b|u|\right))\right)+&2\widetilde A\left(\frac{5b}{\lambda}\widetilde A^{-1}(A\left(b|\nabla u|\right))\right)\bigg]\nonumber\\
	\leq\frac{1}{5}\widetilde A
	\left(\frac{10q(x)}{\lambda}\right)+c_1(1+2F \left(bg_1(x)\right)	+2A\left(b|\nabla u|\right))&:=v_\lambda (x)\,.\nonumber
\end{align}
The function $v_\lambda\in L^1(\Omega)$ because of \eqref{intfinito1} and standard arguments used several times in this proof. Thus, arguing as for \eqref{continuity}, $\abar(x,u_k,\nabla u)-\abar(x,u,\nabla u)\to 0$ in $L^{\widetilde A}(\Omega)$ and 
\begin{equation}\label{ellconv}
	\lim_{k\to +\infty}\int_\Omega (\abar(x,u_k,\nabla u)-\abar(x,u,\nabla u))(\nabla u_k-\nabla u)dx=0\,.
\end{equation}
By definition of weak convergence 
\begin{equation}\label{2}
	\lim_{k\to +\infty}\int_\Omega \abar(x,u,\nabla u)(\nabla u_k-\nabla u)dx=0\,.
\end{equation}
From \eqref{a2}
\begin{align}\label{mag22}
0\leq& \int_\Omega (\abar(x,u_k,\nabla u_k)-\abar(x,u_k,\nabla u))(\nabla u_k-\nabla u)dx\\
 =\int_\Omega (\abar(x,u_k,\nabla u_k)-\abar(x,u,\nabla u))(\nabla u_k-\nabla u)dx&
+\int_\Omega (\abar(x,u,\nabla u)-\abar(x,u_k,\nabla u))(\nabla u_k-\nabla u)dx \,.\nonumber
\end{align}
Passing to the limit in \eqref{mag22} and using \eqref{1}, \eqref{2} and \eqref{ellconv} 
\begin{equation}\label{strongellconv}
	\lim_{k\to +\infty}\int_\Omega (\abar(x,u_k,\nabla u_k)-\abar(x,u_k,\nabla u))(\nabla u_k-\nabla u)dx=0\,.
\end{equation}
Thus the sequence $(\abar(x,u_k,\nabla u_k)-\abar(x,u_k,\nabla u))(\nabla u_k-\nabla u)\to 0$ in $L^1(\Omega)$, and we can find a subsequence, say still $\{u_k\}$, and a set $\Omega_0\subset \Omega$, such that $|\Omega_0|=0$ and (recall that $u_k\to u$ in $L^{F_1}(\Omega)$)
\begin{align}\label{to 0}
	\abar(x,u_k,\nabla u_k)-\abar(x,u_k,\nabla u))(\nabla u_k-\nabla u)\to 0,\ u_k(x)\to u(x) \ \hbox{in}\ \Omega\setminus\Omega_0\,.
\end{align}
 We prove that for every $x\in\Omega\setminus\Omega_0$ there exists $M>0$ such that $|\nabla u_k(x)|\leq M$ for all $k\in N$. We argue by contradiction. Assume that there exists $x\in\Omega\setminus\Omega_0$ such that for every $h>|\nabla u(x)|+1$ there exists $k_h\in N$ such that $|\nabla u_{k_h}(x)|>h$. In particular 
\begin{align}\label{>1}
	|\nabla u_{k_h}(x)-\nabla u(x)|>1\,. \end{align} 
The sequence $\left\{\frac{\nabla u_{k_h}(x)-\nabla u(x)}{|\nabla u_{k_h}(x)-\nabla u(x)|}\right\}$ converges to $\xi\in \rn$, up to a subsequence. We keep the same notation as above for the subsequence and use \eqref{a2} and \eqref{>1} 
\begin{align}\label{zero}
	0\leq  \left(\abar\left(x,u_{k_h}(x),\nabla u(x)+\frac{\nabla u_{k_h}(x)-\nabla u(x)}{|\nabla u_{k_h}(x)-\nabla u(x)|}\right)-\abar(x,u_{k_h}(x),\nabla u(x))\right)\frac{\nabla u_{k_h}(x)-\nabla u(x)}{|\nabla u_{k_h}(x)-\nabla u(x)|}=\\
	\left(\abar\left(x,u_{k_h}(x),\nabla u(x)+\frac{\nabla u_{k_h}(x)-\nabla u(x)}{|\nabla u_{k_h}(x)-\nabla u(x)|}\right)-\abar(x,u_{k_h}(x),\nabla u_{k_h}(x))\right)\frac{\nabla u_{k_h}(x)-\nabla u(x)}{|\nabla u_{k_h}(x)-\nabla u(x)|}\nonumber\\
	+\left(\abar(x,u_{k_h}(x),\nabla u_{k_h}(x))-\abar(x,u_{k_h}(x),\nabla u(x))\right)\frac{\nabla u_{k_h}(x)-\nabla u(x)}{|\nabla u_{k_h}(x)-\nabla u(x)|}\leq \nonumber\\
	\left(\abar(x,u_k(x),\nabla u_k(x))-\abar(x,u_k(x),\nabla u(x))\right)(\nabla u_{k_h}(x)-\nabla u(x))\,.\nonumber
\end{align} 
From \eqref{zero} and \eqref{to 0}
\begin{align}
	\lim _{h\to +\infty}\left(\abar(x,u_{k_h},\nabla u+\frac{\nabla u_{k_h}(x)-\nabla u(x)}{|\nabla u_{k_h}(x)-\nabla u(x)|})-\abar(x,u_{k_h},\nabla u)\right)\frac{\nabla u_{k_h}(x)-\nabla u(x)}{|\nabla u_{k_h}(x)-\nabla u(x)|}=0
\end{align}
and this leads to $\left(\abar(x,u(x),\nabla u(x)+\xi)-\abar(x,u(x),\nabla u(x))\right)\xi=0$, namely $\xi=0$. This contradicts $|\xi|=1$, thus $|\nabla u_{k_h}(x)|\leq M$, for some $M>0$ and we can find a subsequence, say still $\{\nabla u_{k_h}\}$, converging to $\eta \in \rn$. Then from \eqref{to 0} and the convergence of $\{\nabla u_{k_h}\}$ to $\eta \in \rn$
\begin{align}
0=	\lim _{h\to +\infty}\left(\abar(x,u_{k_h}(x),\nabla u_{k_h}(x))-\abar(x,u_{k_h}(x),\nabla u(x))\right)(\nabla u_{k_h}(x)-\nabla u(x))&\\
=\left(\abar(x,u(x),\eta)-\abar(x,u(x),\nabla u(x))\right)(\eta-\nabla u(x))\,.\nonumber
\end{align}
From \eqref{a2} we deduce $\eta=\nabla u(x)$. We have so proved that every subsequence of 
$\{\nabla u_k\}$ has a subsequence converging to $\nabla u(x)$. Thus the whole sequence converges to $\nabla u(x)$ in $\Omega\setminus \Omega _0$.\\ 
Step 2 $|a_i(x,u_k,\nabla u_k)-a_i(x,u,\nabla u)|\rightharpoonup 0$ in $L^{\widetilde{A}}(\Omega)$.\\
Due to step 1 and to the continuity of $\abar(x,\cdot,\cdot)$, we have $\abar(x,u_k,\nabla u_k)\to \abar(x,u,\nabla u)$ a.e. in $\Omega$. $\{u_k\}_{k\in N}$ is bounded in $\w0$ and, from \eqref{abound}, $\{|\abar(x,u_k,\nabla u_k)|\}_{k\in N}$ is bounded in $L^{\widetilde{A}}(\Omega)$. 
Since $\widetilde{A}\in\nabla_2$ near infinity we can apply Corollary \ref{corvit} to obtain 
\begin{equation*}
	|a_i(x,u_k,\nabla u_k)-a_i(x,u,\nabla u)|\rightharpoonup 0 \ \hbox{in}\ L^{\widetilde{A}}(\Omega)\,.
\end{equation*}
Step 3. $\abar(x,u_k,\nabla u_k)\cdot \nabla u_k\to \abar(x,u,\nabla u)\cdot \nabla u$ in $L^{1}(\Omega)$.\\
It holds
\begin{align}\label{sigma}
|\abar(x,u,\nabla u)\cdot\nabla u-\abar(x,u_k,\nabla u_k)\cdot\nabla u_k|=&2 ( \abar(x,u,\nabla u)\cdot\nabla u-\abar(x,u_k,\nabla u_k)\cdot\nabla u_k)^+\\ &-\abar(x,u,\nabla u)\cdot\nabla u+\abar(x,u_k,\nabla u_k)\cdot\nabla u_k	\,.\nonumber
\end{align}
Put $\sigma_k(x)=\abar(x,u,\nabla u)\cdot\nabla u-\abar(x,u_k,\nabla u_k)\cdot\nabla u_k$. One has $\sigma_k^+\leq \abar(x,u,\nabla u)\cdot\nabla u$ and $\sigma_k\to 0$ a.e in $\Omega$. Then 
\begin{align}\label{sigma+}
	\sigma_k(x)^+\to 0\ &\ \hbox{in}\ L^1(\Omega).\end{align} 
From \eqref{sigma}
\begin{align}
0\leq\int_\Omega |\abar(x,u_k,\nabla u_k)\cdot\nabla u_k-\abar(x,u,\nabla u)\cdot\nabla u|dx = &2\int_\Omega \sigma_k^+(x) dx-\int_\Omega \abar(x,u,\nabla u)\cdot\nabla udx\nonumber \\
	+\int_\Omega\abar(x,u_k,\nabla u_k)\left(\nabla u_k-\nabla u\right)dx+&\int_\Omega\abar(x,u_k,\nabla u_k)\nabla u dx \,.\nonumber 
\end{align}
So, using \eqref{sigma+}, \eqref{1} and Step 2
\begin{equation}
	\lim_{k\to +\infty}\int_\Omega |\abar(x,u_k,\nabla u_k)\cdot \nabla u_k- \abar(x,u,\nabla u)\cdot \nabla u|dx=0\,.
\end{equation}
Step 4 $\nabla u_k\to \nabla u$ in $L^A(\Omega)$.\\
From \eqref{a3}, $cA(|\nabla u_k(x)|)\leq \abar(x,u_k,\nabla u_k)\cdot\nabla u_k+dG(d|u_k|)+r(x)$ for all $x\in \Omega$, all $k\in N$. From Step 3 and the assumption $G\ll A_n$, we know that there exists a subsequence of $\{u_k\}$, say still $\{u_k\}$ and a function $g\in L^1(\Omega)$, such that $\abar(x,u_k,\nabla u_k)\cdot\nabla u_k+dG(d|u_k|)\leq g(x)$ a.e. in $\Omega$. Taking into account that $A(c|\nabla u_k(x)|)\to A(c|\nabla u(x)|)$ a.e. in $\Omega$, we deduce that
\begin{equation}\label{finalconv}
\lim_{k\to +\infty}\int_\Omega A(c|\nabla u_k(x)|)dx=\int_\Omega A(c|\nabla u(x)|)dx\,.
\end{equation}
A standard and repeatedly used argument shows that it holds for the whole sequence. Equation \eqref{finalconv} and the $\Delta_2$ condition on $A$, guarantees that $\nabla u_k\to \nabla u$ in $L^A(\Omega)$.\qed
We now construct the truncation of ${\cal S}$ that we use in the proof of our results.
For every $r\in\r$, we set $r^+=\max\{r,0\}$, $r^-=\max\{-r,0\}$.\\
Let $\underline{u},\,\overline{u}\in W^{1,A}(\Omega)$ be such that $(\overline{u})^-,\,\underline{u}^+\in \w0$, and $\underline{u}\leq \overline{u}$ a.e. in $\Omega$. The truncation operator $T:\w0\to \w0$ is defined by
\begin{equation}\label{truncation}
	T(u)=\left\{\begin{array}{cc}
		\overline{u} & \hbox{if}\ u>\overline{u}\\
		u & \hbox{if}\ \underline{u}\leq u\leq \overline{u}\\
		\underline{u} & \hbox{if}\ u<\underline{u}
	\end{array}\right.
\end{equation}
The properties of  $\underline{u},\,\overline{u}$ guarantee that $T(u),\, (u-\overline{u})^+, \,(u-\underline{u})^-\in \w0$. In particular $T$ is well defined.\\
It is known (see \cite {H}, p.20) that
\begin{equation}\label{gradtruncation}
	\nabla T(u)=\left\{\begin{array}{cc}
		\nabla\overline{u}(x) & \hbox{a.e. on the set}\ \{u>\overline{u}\}\\
		\nabla u(x)& \hbox{a.e. on the set}\ \{\underline{u}\leq u\leq \overline{u}\}\\
		\nabla \underline{u}(x) & \hbox{a.e. on the set}\ \{u<\underline{u}\}
	\end{array}\right.
\end{equation}
Given the functions $\underline{u},\,\overline{u}\in W^{1,A}(\Omega)$ as above, and such that $\underline{u},\,\overline{u}\in L^{A_n}(\Omega)$, let us define the operator  ${\cal S}_T: W_0^{1,A}(\Omega)\to(W_0^{1,A}(\Omega))^*$, as
\begin{equation}\label{STdef}
	\langle {\cal S}_{T}u,v\rangle=\int _\Omega \abar(x,Tu,\nabla u)\cdot\nabla v \,dx\,.
\end{equation}
\begin{corollary}\label{STproperties} Assume that $A:[0,+\infty[\to 0,+\infty[$  is a  Young function, $A\in \Delta_2\cap\nabla_2$ near infinity, and that $\abar$ satisfies \eqref{a1}, \eqref{a2} and \eqref{a3bis}. Then the operator ${\cal S}_T$, introduced in \eqref{STdef} is well defined, bounded, continuous and has the $(S_+)$ property.
\end{corollary}
{\bf Proof}. 
The inequality $|T(u(x))|\leq |u(x)|$ for all $x\in \Omega$ guarantees that the operator ${\cal S}_T$ is well defined and bounded. Thanks to Lemma 4.1 in \cite{BaTo1}, the arguments used in Proposition \ref{S1properties} work also for ${\cal S}_{T}$. Thus ${\cal S}_{T}$ is continuous and has the $(S)_+$ property.
\qed

\section{Main results}\label{sec3}
In this Section we state two of the main results of the paper (Theorems \ref{Tnew} and  \ref{Teo1}).\\
First we give the fundamental definitions of weak solution, subsolution and supersolution to \eqref{problem}.
\begin{definition}\label{sol}
	A function  $u\in \w0$ is a weak solution to problem \eqref{problem} if
	\begin{equation*}
		\int_{\Omega}\abar(x,u,\nabla u)\cdot \nabla v dx= \int_{\Omega}f(x,u,\nabla u)vdx\ \hbox{for all $v\in \w0$},
	\end{equation*} 
and $\int_{\Omega}\abar(x,u,\nabla u)\cdot \nabla v dx\in \r$ for all $v\in \w0$.
\end{definition}
\begin{definition}\label{supersol}
We say that $\overline u\in W^{1,A}(\Omega)$ is a supersolution to \eqref{problem} if  $(\overline{u})^{-} \in \w0$,
\begin{equation*}\label{supersolution}
	+\infty>\int_{\Omega}\abar(x,\overline u,\nabla \overline u)\cdot\nabla v dx\geq \int_{\Omega}f(x,\overline u,\nabla \overline u)vdx>-\infty
\end{equation*}
for all $v\in \w0$, $v\geq 0$ a.e. in $\Omega$.
\end{definition}
\begin{definition}\label{subsol}
	We say that $\underline u\in W^{1,A}(\Omega)$ is a subsolution to \eqref{problem} if $\underline{u}^+\in \w0$ and
\begin{equation*}\label{subsolution}
-\infty<	\int_{\Omega}\abar(x,\underline u,\nabla\underline  u)\cdot\nabla v dx\leq \int_{\Omega}f(x,\underline u,\nabla \underline u)vdx<+\infty
\end{equation*}
for all $v\in \w0$, $v\geq 0$ a.e. in $\Omega$.
\end{definition}
\begin{theorem}\label{Tnew}
Let $\Omega$ be an open set in $\rn$, with $n\geq 2$, such that $|\Omega|<\infty$. Let $A\in C^1([0,+\infty))$ be a Young function, $A\in \Delta_2\cap\nabla_2$ near infinity. Assume also that $A$ satisfies \eqref{conv0} and \eqref{div}, or \eqref{conv}. Let $\underline{u}$ and $\overline{u}$ be a subsolution and a supersolution of problem (\ref{problem}), respectively, with $\underline{u}\leq\overline{u}$ a.e. in $\Omega$, and $\underline{u},\overline{u}\in L^{A_n}(\Omega)$. Assume that the function $\abar$ satisfies \eqref{a1}, \eqref{a2}, \eqref{a3}. Let $f:\Omega\times \mathbb{R}\times \mathbb{R}^n\to \mathbb{R}$ be a Carath\'{e}odory function fulfilling
\begin{equation}\label{growth f'}
|f(x,s,\xi)|\leq \sigma (x)+\overline \gamma \widetilde E^{-1}(A(|\xi|))\ \hbox{  for a.e.}\ x\in \Omega,\ \hbox{all}\ s\in [\underline u(x),\overline u(x)],\ \hbox{all}\ \xi\in \rn\,,
\end{equation}
where $\sigma\in L^{\widetilde A_n}(\Omega)$, $\overline \gamma>0$ and $E:[0,+\infty[\to [0,+\infty[$ is a Young function, $E\ll A_n$ near infinity.\\
Then problem $(P)$ has a solution $u\in \w0$ such that $\underline{u}\leq u\leq \overline{u}$ a.e. in $\Omega$.
\end{theorem}
To prove Theorem \ref{Tnew}, we perturb problem (\ref{problem}) and formulate an auxiliary one. Let
$\Pi:W_0^{1,A}(\Omega)\to( W_0^{1,A}(\Omega))^*$, given by
\begin{equation}\label{Pi}
	\Pi(u)(v)=\int_\Omega\pi(x,u(x))v(x)dx\,,\ \hbox{for}\ u,v\in \w0,
\end{equation}
where 
\begin{equation*}
	\pi(x,s)=\left\{\begin{array}{cc}
		\widetilde E^{-1}(E(s-\overline u(x) )& \hbox{if}\ s>\overline u (x)\\
		0& \hbox{if}\ \underline u(x)\leq s\leq\overline u (x)\\
		-\widetilde E^{-1}(E(\underline u(x)-s)) & \hbox{if}\ s<\underline u(x).
	\end{array}\right.
\end{equation*}
Let ${\cal N}_f\circ T:\w0\to(\w0)^*$ be the operator defined as
\begin{equation*}
\langle{\cal N}_f\circ T(u),v\rangle=\int_\Omega f(x,Tu,\nabla Tu)v(x) dx\,,\ \hbox{for}\ u,v\in \w0,
\end{equation*}
Given $\mu >0$, we consider the following problem

\begin{equation} \label{Plambda}
	\begin{cases}
		- {\rm div}(\abar(x,Tu,\nabla u))+\mu\Pi
		(u)=N_f(Tu)&\hbox{in}\; \Omega,\\
		u=0&\hbox{on}\; \partial \Omega.
	\end{cases}
\end{equation}
The result below guarantees that problem \eqref{Plambda} has a solution, provided the parameter $\mu>0$ is sufficiently large.

\begin{theorem}\label{T1}

Under the same assumtions of Theorem \ref{Tnew}, there exists $\mu_0>0$ such that \eqref{Plambda} admits a solution whenever $\mu\geq \mu_0$.
\end{theorem}
{\bf Proof.} For all $\mu>0$ consider the operator ${\cal A}_\mu :W_0^{1,A}(\Omega)\to (W_0^{1,A}(\Omega))^*$, defined by
\begin{align*}
	\langle {\cal A}_\mu (u),v\rangle=&\int_\Omega \abar(x,Tu,\nabla u)\cdot\nabla v \,dx+\mu
	\int_\Omega \pi(x,u)v\,dx-\int_\Omega f(x,Tu,\nabla Tu)v\,dx=\\
	&\langle {\cal S
	}_Tu+\mu\Pi u-{\cal N}_f\circ T(u),v \rangle\qquad\qquad \hbox{for }\ u,v\in W_0^{1,A}(\Omega)\,.
\end{align*}
We prove that ${\cal A}_\mu$ is well defined, bounded, pseudomonotone and there is $\mu_0>0$ such that ${\cal A}_\mu$ is coercive for all $\mu>\mu_0$.\\
Due to Corollary \ref{STproperties}, Propositions 4.3 and 4.5 of \cite{BaTo1}, ${\cal A}_\mu$ is well defined, bounded and continuous. To prove that it is pseudomonotone, we take $u\in \w0$, and a sequence  $\{u_k\}\subset \w0$ such that
\begin{equation*}
	u_k\rightharpoonup u \quad\hbox{in $\w0$},\ \hbox{and}\ \limsup_{k\to\infty}\,\langle{\cal A}_\mu(u_k), u_k-u \rangle\leq 0\,.
\end{equation*}
Equations (4.6) and (4.17) of \cite{BaTo1} allow to write
\begin{equation*}
	\limsup_{k\to \infty} \, \langle {\cal S}_{T}(u_k), u_k-u\rangle\leq 0\,.
\end{equation*}
Thus $u_k\to u$ in $\w0$ (see Corollary \ref{STproperties}), and consequently 

\begin{equation*}
	\lim_{k\to \infty}\|{\cal A}_\mu(u_k)-{\cal A}_\mu(u)\|_{(\w0)^*}=0\,,
\end{equation*}
so $\langle{\cal A}_\mu(u_k),u_k\rangle\to \langle {\cal A}_\mu(u),u\rangle$, $\langle{\cal A}_\mu(u_k),v\rangle\to \langle {\cal A}_\mu(u),v\rangle$ for all $v\in \w0$ and ${\cal A}_\mu$ is a pseudomonotone operator.\\
It remains to prove that ${\cal A}_\mu$ is coercive for some $\mu>0$. Arguing like for equation (5.3) of \cite{BaTo1} we can find a constant $c_1>0$ such that
\begin{equation}\label{c1}
	\overline \gamma\int_\Omega \widetilde E^{-1}(A(|\nabla Tu|))|u|dx\leq c_1+2\overline \gamma\int_\Omega E\left(\frac{|u|}{2}\right)dx +\frac{c}{2}\int_\Omega A(|\nabla u|)dx\ \hbox{for all}\ u\in\w0\,.
\end{equation}
Here $c$ is that of \eqref{a3}. From \eqref{growth f'}, \eqref{holderyoung}, \eqref{B-Wbis}, and \eqref{c1}
\begin{eqnarray}\label{fpercoerc}
	\left|\int_\Omega f(x,Tu,\nabla Tu)udx\right|\leq\left|\int_\Omega\sigma (x)u(x)dx\right|+\overline \gamma\left|\int_\Omega \widetilde E^{-1}(A(|\nabla Tu|))udx\right|\leq \nonumber \\
	\leq 2C\|\sigma\|_{L^{\widetilde B}(\Omega)} \|u\|_{W_0^{1,A}(\Omega)}+\frac{c}{2}\int_{\Omega}A(|\nabla u|)dx+2\overline \gamma\int_\Omega E\left(\frac{|u|}{2}\right)dx+c_1\,.
\end{eqnarray}
Let $v(x)=\max\{|\underline u(x)|, |\overline u(x)|\}$ Then $|T(u)(x)|\leq v(x)$ for all $x\in \Omega$, and $\int_\Omega G(dv(x))dx<+\infty$. From \eqref{a3} 
\begin{align}\label{pAc}
\int_\Omega \abar(x,Tu,\nabla u)\cdot\nabla u\,dx&\geq c\int_\Omega A(|\nabla u|) dx-d\int_\Omega G(d|T(u)|)dx-\int_\Omega r(x)dx\\
&\geq c\int_\Omega A(|\nabla u|)dx-c_2\ \ \hbox{for any}\ u\in \w0\,.\nonumber 
\end{align}
So, choosing $\mu >a$, $u\in \w0$ with $\|u\|_{\w0} >> 1$, and using \eqref{fpercoerc}, \eqref{pAc}, Lemma 4.6 of \cite{BaTo1},  and \eqref{min per coerc}
\begin{align*}
\frac{\langle{\cal A}_\mu(u),u\rangle}{\|u\|_{\w0}} \geq&
\frac{\frac{c}{2}\int_\Omega A(|\nabla u|)dx -2C\|\sigma\|_{L^{\widetilde B}(\Omega)} \|u\|_{\w0}+2(\mu-\overline \gamma)\int_\Omega E\left(\frac{|u|}{2}\right)dx -c_3}{\|u\|_{\w0}}\\
&\geq
\frac{\frac{ck_1^{i^\infty_A}}{2}\|u\|_{\w0}^{i_A^\infty} -2C\|\sigma\|_{L^{\widetilde B}(\Omega)} \|u\|_{\w0}+2(\mu-\overline \gamma)\int_\Omega E\left(\frac{|u|}{2}\right)dx -c_4}{\|u\|_{\w0}}\,.\nonumber
\end{align*}
Thus
\begin{equation*}
	\lim_{\|u\|\to +\infty}\frac{\langle {\cal A}_\mu(u),u\rangle}{\|u\|_{\w0}}=+\infty\,.
\end{equation*}
Theorem \ref{theorem on pseudomonotone operator} guarantees that there exists $u\in W^{1,A}_0(\Omega)$ such that ${\cal A}_\mu(u)\equiv 0$. Thus
\begin{equation}\label{solution}
	\int_{\Omega}\abar(x,Tu,\nabla u)\cdot  \nabla vdx+\mu \int_{\Omega} \pi(x,u(x))v(x)dx- \int_{\Omega}f(x,Tu,\nabla Tu)vdx=0\,.
\end{equation} for all $v\in W_0^{1,A}(\Omega)$.\qed
\begin{remark}
	The proof above works also if we weaken \eqref{a1}$\ldots$\eqref{a3bis}, requiring them to hold for $s\in [\underline u(x), \overline u(x)]$ rather than for all $s\in \r$. This is because in the proof we consider only the truncated function $\abar(x,Tu,\nabla u)$. 
\end{remark}
{\bf Prooof of Theorem \ref{Tnew}.} By Theorem \ref{T1} there exists a solution $u\in \w0$ of the truncated auxiliary problem \eqref{Plambda} provided $\mu>0$ is sufficiently large. Let us fix such a $\mu$ and $u$.\\
Via the same comparison arguments of the proof of Theorem 3.6 of \cite{BaTo1} we can prove that the solution of \eqref{Plambda} has the enclosure property $u\in [\underline{u},\overline{u}]$. Thus, it follows from \eqref{truncation} and \eqref{Pi} that $Tu=u$ and $\Pi(u)=0$. Consequently, $u$ is a solution of\eqref{problem}.\qed
The proof of the Corollary below follows the same lines as that of Corollary 5.2 
of \cite{BaTo1}.
\begin{corollary}\label{Cor1}
Let $\Omega$ be an open set in $\rn$, with $n\geq 2$, such that $|\Omega|<\infty$. Let $A\in C^1([0,+\infty))$ be a Young function, $A\in \Delta_2\cap\nabla_2$ near infinity. Assume also that $A$ satisfies \eqref{conv0} and \eqref{div}, or \eqref{conv}. Let $\underline{u}$ and $\overline{u}$ be a subsolution and a supersolution of problem (\ref{problem}), respectively, with $\underline{u}\leq\overline{u}$ a.e. in $\Omega$, $\underline{u},\overline{u}\in \w0$, and such that the Carath\'{e}odory function $f:\Omega\times \mathbb{R}\times \mathbb{R}^n\to \mathbb{R}$ fulfills
\begin{equation}\label{growth fcor1}
		|f(x,s,\xi)|\leq \rho(x)+g(|s|)+\overline{\gamma}\widetilde E^{-1}(A(|\xi|))\ \hbox{a.e.}\ x\in \Omega,\ \hbox{all}\ s\in [\underline{u}(x),\overline{u}(x)], \ \xi\in\rn\,,
	\end{equation}
where $\rho\in L^{\widetilde A_n}(\Omega)$, $\overline{\gamma},\,E,$ are as in Theorem \ref{Tnew}, and $g:[0,+\infty[\to [0,+\infty[$ is a nondecreasing function such that $g(0)=0$ and there exist $s_0,\,h>0$ such that $g(|s|)|s|\leq A_n(h|s|)$ for all $|s|\geq s_0$.\\
\noindent Then problem $(P)$ possesses a nontrivial solution $u\in \w0$.
\end{corollary}

We consider now a special instance of \eqref{problem}, in which $\abar$ does not depend on $s$ and has a potential with respect to $\xi$. So, let $\Omega\subset\rn$ be a set of finite measure and let $A,B$ be two Young functions such that $A\in \Delta_2\cap \nabla_2$ near infinity, $B\in\nabla_2$ near zero, and $A\circ B^{-1}$ is a Young function too. We assume that
$\abar:\Omega\times\rn\to\rn$, $\abar=(a_1,\ldots a_n)$, is such that each $a_i(x,\xi)$ is a Carath\'{e}odory function, and
\begin{align}\label{a_1'}
	&|\abar(x,\xi)|\leq q(x)\widetilde B^{-1}(B(b|\xi|))+ b\widetilde A^{-1}(A(b|\xi|))\\ 
	&\hbox{for some}\ q\in L^{\widetilde{A\circ B^{-1}}}(\Omega),\ \hbox{some} \ b>0,\, \hbox{for a.e.} \, x\in\Omega,\, \hbox{all}\ \xi\in\rn\,,\nonumber
\end{align}
\begin{align}\label{a_2'}
	\sum_{i=1}^n\left(a_i(x,\xi)-a_i(x,\xi ')\right)\cdot (\xi_i-\xi_i')>0 \quad \hbox{for a. e.}\ x\in \Omega,\ \hbox{all}\ \xi,\xi'\in\rn\,,\ \xi\neq\xi'\,,
\end{align}
\begin{align}\label{a_3'}
\sum_{i=1}^n a_i(x,\xi)\cdot \xi_i\geq cA(c|\xi|)\quad \hbox{for some}\ c>0,\ \hbox{for a.e.}\ x\in\Omega,\ \hbox{all}\ \xi\in\rn\,.
\end{align}
Furthermore, we assume that there exists a measurable function $\Phi:\Omega\times\rn\to\r$, even with respect to $\xi\in\rn$ and such that 
\begin{align}\label{a4}	
	\Phi_\xi(x,\xi)=\abar(x,\xi)\ \hbox{for all} \ (x,\xi)\in \Omega\times \rn,&\  \Phi(x,0)=0\ \hbox{for all} \ x\in\Omega.
\end{align}
Since $A\circ B^{-1}$ is a Young function, it follows that $A$ dominates $B$ near infinity and $B$ dominates $A$ near zero. Thus
\begin{align}\label{ABdomin}
	A(t)\leq B(\overline kt)\ \hbox{for}\ 0\leq t\leq \overline t \quad \hbox{and} &\quad B(t)\leq A(\tilde kt)\ \hbox{for}\ t\geq \tilde t>0\,.
\end{align} 
Also, if $u\in L^A(\Omega)$ then $B(k|u|)\in L^{A\circ B^{-1}}(\Omega)$ for all $k>0$.\\
Condition \eqref{a_2'} ensures that $\Phi(x,\cdot)$ is convex for every $x\in \Omega$.
From \eqref{a_1'} and \eqref{a_3'} there exist $k_4,k_5>0$ such that 
\begin{equation}\label{equiv}
	k_4A\left(k_4|\xi|\right)\leq\Phi(x,\xi) \leq  2q(x)B(b|\xi|)+ k_5A(k_5|\xi|) \ \hbox{for all}\ (x,\xi)\in\rn\,.
\end{equation}
Now, problem \eqref{problem} reads as 
\begin{equation}\label{problemvar}
	\begin{cases}
		- {\rm div}(\Phi_\xi(x,\nabla u)) =f(x,u, \nabla u) & {\rm in}\,\,\, \Omega \\
		u  =0 & {\rm on}\,\,\,
		\partial \Omega \,.
	\end{cases}
\end{equation} 
For functions $f$ satisfying suitable growth conditions we can construct a sub  or a supersolution for problem \eqref{problemvar}, via variational methods. The $\Delta_2$ and $\nabla_2$ conditions play a crucial role here.
\begin{theorem}\label{Teo1}
Let $\Omega$ be an open set in $\rn$, with $n\geq 2$, such that $|\Omega|<\infty$. Let $A\in C^1([0,+\infty))$ be a Young function, $A\in \Delta_2\cap \nabla_2$ at infinity. Assume also that $A$ satisfies \eqref{conv0} and \eqref{div}, or \eqref{conv}. Let $\abar: \Omega\times \rn \to \rn$ and $\Phi: \Omega\times \rn \to \r$ be two Carath\'{e}odory functions satisfying  \eqref{a_1'}$\ldots$\eqref{a4}. Let $f:\Omega\times \mathbb{R}\times \mathbb{R}^n\to \mathbb{R}$ be a Carath\'{e}odory function fulfilling
\begin{align}\label{growth f1}
-\rho_1(x)-g_1(|s|)\leq f(x,s,\xi)\leq&\rho_2(x)+g_2(|s|)+\overline \gamma\widetilde E^{-1}(A(|\xi|))\ \hbox{for a.e.}\, x\in \Omega,\ \hbox{all}\ s\leq 0,\\
\ \hbox{all}\ \xi\in\rn,\ f(x,0,0)\leq 0\ \hbox{in}\ \Omega&\ \hbox{and}\ f(x,0,0)<0\ \hbox{on a set of positive measure},\nonumber 		
	\end{align}
or
\begin{align}\label{growth f2}
-\rho_2(x)-g_2(|s|)-\overline \gamma\widetilde E^{-1}(A(|\xi|))&\leq f(x,s,\xi)\leq \rho_1(x)+g_1(|s|)\ \hbox{for a.e.}\ x\in \Omega,\ \hbox{all}\ s\geq 0,\\
  \ \hbox{all}\ \xi\in \rn,\ f(x,0,0)\geq 0\ \hbox{in}\ \Omega,&\ \hbox{and}\ f(x,0,0)>0\ \hbox{on a set of positive measure},\nonumber 
\end{align}
where $\overline \gamma>0$, $E$ is a Young function, $E\ll A_n$ near infinity,
	$\rho_1,\rho_2:\Omega\to [0,+\infty[$ are two measurable functions, $\rho_i\in L^{\widetilde A_n}(\Omega),\ i=1,2,$, $g_1,g_2:[0,+\infty[\to [0,+\infty[$ are two non decreasing functions such that $g_1(0)=g_2(0)=0$
and there exist $s_0>0,\,h_0\in \left ]0,\tau\omega_n^{\frac{1}{n}}|\Omega|^{-\frac{1}{n}}\right[$, $h_1>0$ such that
\begin{align}\label{g12}
	g_1(|s|)|s|\leq A(h_0|s|)\ \hbox{and}\  g_2(|s|)|s|\leq A_n(h_1|s|)\ \hbox{for all}\ |s|\geq s_0\,.
\end{align}	
Here $\omega_n$ is the measure of the unit ball in $\rn$, $\tau=\min\{1,k_4^2\}$ where $k_4$ is that of \eqref{equiv}.\\
Then problem \eqref{problemvar} possesses a nontrivial constant sign solution $u\in \w0$. 
\end{theorem}
{\bf Proof.} Suppose that \eqref{growth f1} is in force. We construct a subsolution $\underline u\leq 0$ a.e., $\underline u\not\equiv 0$, and show that $\overline u\equiv 0$ is a supersolution but not a solution to \eqref{problemvar}. Then, we show that $f$ satisfies \eqref{growth fcor1}.\\
Put $G_1(t)=\int_0^t g_1(\tau)d\tau,\ t\geq 0$ and consider the functional $J:\w0\to \r$, defined as
\begin{equation*}\label{J}
	J(u)=\int_\Omega \left(\Phi(x,\nabla u)+\rho_1(x)u-G_1(|u|)\right)dx\quad \hbox{for}\ u\in \w0.
\end{equation*}
We prove that $J$ is well defined, weakly lower semicontinuous, coercive and 
\begin{equation}\label{J'}
	J'(u)v= \int_{\Omega}\Phi_\xi(x,\nabla u)\nabla v dx+\int_{\Omega}\rho_1(x)v(x)dx-\int_{\Omega}g_1(|u|)sign\, u\, v(x)dx
\end{equation}
for all $u,v\in\w0$. We examine separately the three integrals.\\
Due to \eqref{equiv}, the fact that $A\in\Delta_2$ at infinity, and the convexity of $\Phi(x,\cdot)$, for all $x\in\Omega $, the functional $u\mapsto \int_\Omega \Phi(x,\nabla u)dx$
is well defined in $\w0$, convex. We briefly sketch the proof of its regularity, because it makes use of standard arguments like the Lebesgue Theorem, and the properties of Young's functions. Let $u,v\in\w0$. For all $x\in\Omega$, all $t>0$, $t<<1$, there exists $\mu_{t,x}\in (0,1)$ such that 
\begin{align}
\left|\frac{\Phi(x,\nabla u+t\nabla v)-\Phi(x,\nabla u) }{t}\right|=&\left|\Phi_\xi(x,\nabla u+\mu_{t,x}t\nabla v)\nabla v\right|\\
\leq q(x)\widetilde B^{-1}(B(b|\nabla u+\mu_{t,x}t\nabla v|))|\nabla v|+& b\widetilde A^{-1}(A(b|\nabla u+\mu_{t,x}t\nabla v|))|\nabla v|\nonumber \\
\leq 2q(x)\frac{B(b(|\nabla u|+|\nabla v|))}{b(|\nabla u|+|\nabla v|)}|\nabla v|+&2\frac{A(b(|\nabla u|+|\nabla v|))}{|\nabla u|+|\nabla v|}|\nabla v|\nonumber \\
\leq \frac{2q(x)}{b}B(b(|\nabla u|+|\nabla v|))+&2bA(b(|\nabla u|+|\nabla v|))\,.\nonumber
\end{align}	
We used \eqref{a_1'}, the monotonicity of $\widetilde A^{-1}\circ A$ and $\widetilde B^{-1}\circ B$, and \eqref{A1}. Now, the 
condition $A\in\Delta_2$ near infinity, \eqref{ABdomin} and the 
Lebesgue Theorem, allow to prove that the functional $u\mapsto \int_\Omega \Phi(x,\nabla u)dx$ is $C^1$.\\
Thus the weak lower semicontinuity of $J$ and equation \eqref{J'} follow.\\
To prove the coercivity of $J$ we need the following inequality, that can be found in \cite[Proposition 3.2]{BaCi1},
\begin{equation*}\label{B-W}
	\int_{\Omega}A(|u|)dx \leq \int_{\Omega}A(\omega_n^{-\frac{1}{n}}|\Omega|^{\frac{1}{n}}|\nabla u|)dx\ \hbox{for all}\ u\in\w0\,.
\end{equation*}
Let $\ep>0$ be such that $h_0\omega_n^{-\frac{1}{n}}|\Omega|^{\frac{1}{n}}<\tau-\ep$. Using \eqref{g12} and the inequality above
\begin{align}\label{Gc}
	\int_\Omega G_1(|u|)dx\leq &\int_{\{|u|\leq s_0\}} G_1(s_0)dx+\int_{\{|u|>s_0\}} A(h_0|u|)dx \\
	\leq G_1(s_0)|\Omega|+\int_{\Omega}A((\tau-\ep)|\nabla u|)dx\leq&
	 G_1(s_0)|\Omega|+\sqrt{\tau-\ep}\int_{\Omega}A(k_4|\nabla u|)dx\ \hbox{for all}\ u\in\w0\,.\nonumber
\end{align}
Take now $u\in\w0$, $\|u\|_{\w0}>1$ and use \eqref{equiv}, \eqref{holderyoung}, \eqref{Gc} and \eqref{min per coerc} 
\begin{align*}
\frac{J(u)}{\|u\|_{\w0}}&\geq \frac{\left(k_4-\sqrt{\tau-\ep}\right)\int_\Omega A(k_4|\nabla u|)dx -c_4\|u\|_{\w0}-G_1(s_0)|\Omega|}{\|u\|_{\w0}}\\
&\geq \left(k_4-\sqrt{\tau-\ep}\right)(k_4)^{i^\infty_A}\|u\|_{\w0}^{i_A^\infty-1}-c_4-\frac{k_3+G_1(s_0)|\Omega|}{\|u\|_{\w0}}\,.
\end{align*}
This proves that $J$ is coercive. Thus it has a global minimum. Let $\underline u$ be a global minimum point for $J$. We prove that $\underline u\not\equiv 0$. To this end consider a function $v\in C_0^1(\Omega)$, such that $b|\nabla v(x)|\leq \overline t$ and $k_5\overline k|\nabla v(x)|\leq \overline t$ for all $x\in\Omega$. Also, $v\leq 0$ and $\rho_1(x)v(x)\not\equiv 0$ in $\Omega$. The inequality $\frac{B(t_1)}{B(t_0)}>\left(\frac{t_1}{t_0}\right)^{k_B}$ holds for $0<t_0<t_1<\overline t$, and some $k_B>1$, by virtue of the $\nabla_2$ condition near zero. Then, choosing once $t_1=b|\nabla v|$, $t_0=bt|\nabla v|$, and  secondly $t_1=k_5\overline k|\nabla v|$, $t_0=tt_1$, with $t<1$, and taking into account \eqref{equiv}

\begin{eqnarray*}
J(tv)\leq 2t^{k_B}\int_\Omega q(x) B(b|\nabla v|)dx+k_5 t^{k_B}\int_\Omega B(k_5\overline k|\nabla v|)dx+t\int_\Omega \rho_1(x)v dx<0\ \ \hbox{for}\ t <<1\,,
\end{eqnarray*}
and this proves that $J(\underline u)<0$. Using $J(-|\underline u|)\geq J(\underline u)$ and the fact that $\Phi(x,\cdot)$ is even, we obtain $\underline u\leq 0$ a.e. in $\Omega$.\\
Now we prove that $\underline u$ is a subsolution and $u\equiv 0$ is a supersolution but not a solution to \eqref{problem}. Note that 
\begin{equation}\label{J'equiv0}
J'(\underline u)(v)=\int_{\Omega}\Phi_\xi(x,\nabla\underline u)\nabla v dx+\int_{\Omega}(\rho_1(x)+g_1(|\underline u(x)|))vdx
\equiv 0,\ \hbox{for all}\ v\in \w0\,.
\end{equation}
Acting with any $v\in \w0$, $v\geq 0$, in \eqref{J'equiv0} and using \eqref{growth f1}
\begin{equation*}\label{}
	\int_{\Omega}\Phi_\xi(x,\nabla\underline u)\nabla v dx- \int_{\Omega}f(x,\underline u,\nabla\underline u)vdx\leq 0\,,
\end{equation*}
that is $\underline u$ is a subsolution to \eqref{problemvar}. Using \eqref{growth f1} and choosing $v\in \w0$, $v\geq 0$
\begin{equation*}\label{0super}
	0-\int_{\Omega}f(x,0,0)vdx\geq 0\,,
\end{equation*}
thus $u\equiv 0$ is a supersolution to \eqref{problem} and the assumption on $f(x,0,0)$ guarantees that it is not a solution.\\
We put $\rho(x)=\max\{\rho_i(x),\ i=1,2\}$, $g(|s|)=\max\{g_i(|s|),\ i=1,2\}$ and use \eqref{growth f1}
\begin{equation*}\label{}
	|f(x,s,\xi)|\leq\rho(x)+g(|s|)+\overline \gamma\widetilde E^{-1}(A(|\xi|))\ \ \hbox{for}\ x\in \Omega,\ s\in [\underline u(x),0],\ \xi\in \rn \,.
\end{equation*}
Then $f$ satisfies \eqref{growth fcor1} and from Corollary \ref{Cor1} problem \eqref{problemvar} has a nontrivial solution $u\in\w0$ and $u\in [\underline u,0]$.\\
When \eqref{growth f2} is in force we consider $f_1(x,s,\xi)=-f(x,-s,-\xi)$. Then, by virtue of the proof above, problem 
\begin{equation}\label{problemf1}
	\begin{cases}
		- {\rm div}\left(\Phi_\xi(x,\nabla v)\right) =f_1(x,v,\nabla v) & {\rm in}\,\,\, \Omega \\
		v  =0 & {\rm on}\,\,\,
		\partial \Omega \,,
	\end{cases}
\end{equation}
has a nontrivial solution $v\in\w0$, $v\leq 0$ a.e. in $\Omega$. Then the function $u=-v$ is a nontrivial solution to \eqref{problem} and $u\geq 0$ a.e. in $\Omega$. \qed

\begin{remark}\label{rembound}
From \cite[Theorem 3.1]{BaCiMa}, when $\rho_1,\ \rho_2\in L^{M,\infty}(\Omega)$, for a suitable Young function $M$ (see equation $(3.9)$ in \cite{BaCiMa}), then the solution $u$ is essentially bounded.
\end{remark}
\par

\section{Regularity results}\label{sec4}
In this section we give some existence and regularity results, Theorems \ref{reg} and \ref{main1r}. We strenghten the hypotheses on $\Omega$ and on $A$ (see \eqref{dg}), in order to apply regularity theory (see Proposition \ref{Lib}).\\ 
The proof of the existence is based on sub and supersolution methods, while the main tool for the regularity is Theorem 1.7 of \cite{Li1} (see also the remark after that result and \cite{Li}), that we recall below.
\begin{proposition}(see \cite[Theorem 1.7]{Li1})\label{Lib}
	Let $\Omega$ be a bounded domain in $R^n$ with a $C^{1,\alpha}$ boundary, for some $0<\alpha\leq 1$. Let $A$ be a Young function satisfying
\begin{align}\label{dg}
A\in C^2(]0,+\infty[),&\ \hbox{and there exist two positive constants}\	
\delta,\,g_0>0\\
\hbox{such that}\ &\ \delta\leq \frac{tA''(t)}{A'(t)}\leq g_0\quad \hbox{for}\ t>0\,.\nonumber
	\end{align}  
	Let $\abar:\Omega\times \r\times \rn\to \rn $ be a vector valued function, with Carath\'{e}odory components, $a_i$, $i=1,\ldots,n$. Consider the problem
	$$-div(\abar(x,u,\nabla u))=f(x,u,\nabla u)\quad \hbox{in}\ \Omega\,.$$
	Suppose $\abar$ and $f$ satisfy the structure conditions (here $a_{ij}(x,s,\eta)=\frac{\partial a_i}{\partial\eta_j}$) 
	\begin{align}\label{a}
	\sum_{i,j=1}^na_{ij}(x,s,\eta)\xi_i\xi_j\geq  \frac{A'(|\eta|)}{|\eta|}|\xi|^2
	\end{align}
	\begin{align}\label{b}
\sum_{i,j=1}^n|a_{ij}(x,s,\xi)|\leq  \Lambda\frac{A'(|\xi|)}{|\xi|}
\end{align}

	\begin{align}\label{c}
	|\abar(x,s,\xi)-\abar(y,w,\xi)|\leq \Lambda_1(1+A'(|\xi|)(|x-y|^\alpha+|s-w|^\alpha)\end{align}
\begin{align}\label{d}
		|f(x,s,\xi)|\leq \Lambda_1(1+A'(|\xi|)|\xi|),
	\end{align}
	for some positive constants $\Lambda$, $\Lambda_1$, $M_0$, for all $x$ and $y\in \Omega$, for all $s,w\in [-M_0,M_0]$ and for all $\xi \in \mathbb{R}^n$.	
	Then, any solution $u\in W^{1,A}(\Omega)$, with $|u|\leq M_0$ in $\Omega$, is  $C^{1,\beta}(\overline{\Omega})$ for some positive $\beta$.
\end{proposition}
We point out that \eqref{dg} guarantees that $A'(0)=0$ and $A\in \nabla_2\cap\Delta_2$ globally. 

\begin{lemma}\label{LemLib} Let $A$ be a Young function satisfying \eqref{dg}. If $\abar:\Omega\times \r\times \rn\to \rn $ has Carath\'{e}odory components, $a_i(x,s,0)\equiv 0$, for a.e. $x\in\Omega$, all $s\in \r$, all $i=1,\ldots,n$, and satisfies \eqref{a} and \eqref{b} for all $s\in\r$, then \eqref{a1}, \eqref{a2} and \eqref{a3bis} hold with $q(x)\equiv r(x)\equiv 0$, $F\equiv 0$, and $d=0$.
	
\end{lemma}
{\bf Proof}. Let $H(t)=a_i(x,s,t\xi),\ t>0$. Using \eqref{b}, $a_i(x,s,0)\equiv 0$, and \eqref{dg}
\begin{equation}
|H(1)|\leq 	\int_0^1|\nabla a_i(x,s,t\xi)||\xi|dt\leq \Lambda \int_0^1 \frac{A'(t|\xi|)}{t|\xi|}|\xi|dt\leq \frac{\Lambda}{\delta} A'(|\xi|).
\end{equation}
Then, thanks to equation (6.23) in \cite{BaCi1}
\begin{equation}
	|\abar (x,s,\xi)|\leq \frac{\Lambda \sqrt{n}}{\delta}A'(|\xi|)\leq \frac{\Lambda \sqrt{n}}{\delta}\widetilde  A^{-1}(A(2|\xi|)\leq b\widetilde A^{-1}(A(b|\xi|),\ \hbox{for}\ b=\max\left\{\frac{\Lambda \sqrt{n}}{\delta},2\right\}\,,
\end{equation}
and \eqref{a1} is proved. \\
For any $\xi,\xi'\in\rn$, write $\xi=\xi'+\eta$, $\eta\neq 0_{\rn}$. Define
$H(t)=\sum_{i=1}^n\left(a_{i}(x,s,\xi'+t\eta)-a_{i}(x,s,\xi')\right)\eta_i$, for $t\in[0,1]$. From \eqref{a}, $H'(t)\geq \frac{A'(|\xi'+t\eta|)}{|\xi'+t\eta|}\cdot |\eta|^2$ and this guarantees that $H(1)>H(0)=0$, namely \eqref{a2}.\\
Choose now $\xi'=0$ in the function $H$ just used.  
For $\eta\neq 0_{\rn}$, \eqref{a} and \eqref{dg} give 
$$H(1)\geq \int_0^1 \frac{A'(t|\eta|)}{t|\eta|}|\eta|^2dt\geq\int_0^1 \frac{t^{g_0}A'(|\eta|)|\eta|}{t}dt\geq \frac{1}{g_0}A\left(|\eta|\right)\,.$$
If $\eta=0_{\rn}$ then $H(1)=0$, and \eqref{a3bis} is proved.
\qed 
Following the spirit of Example \ref{ex1}, we present a function $\abar$ satisfying \eqref{a}, \eqref{b} and \eqref{c}. This general function will be used in Example \ref{ex3}, but we prefere to introduce this function here, to emphasize its general structure.
\begin{example}\label{abarforreg}
	Let $p>1$, $q\in \r$ and $p+q-1>0$. Consider the Young function $A:[0,+\infty[\to [0,+\infty[$ complying with 
	\begin{align}\label{A'ex}
		A'(t)=t^{p-1}\lg^{q}(1+t)\quad \ \hbox{for}\ t\geq 0\,.
	\end{align}
Then  $A$ satisfies \eqref{dg} because
	\begin{align}\label{A''ex}
		0<\min\{p+q-1, p-1\}\leq\frac{tA''(t)}{A'(t)}&
		\leq \max\{p+q-1, p-1\}<\infty\ \hbox{for}\ t> 0\,.
	\end{align}
	Let us define 
	\begin{align}\label{areg}
		\abar(x,s,\xi)=(\|x\|^\gamma|s|^\delta+1)|\xi|^{p-2}\lg^{q}(1+|\xi|)\xi\ &
		\ \hbox{for}\ (x,s,\xi)\in \Omega\times\r\times\rn\,.
	\end{align}
	Here $\Omega$ is a bounded domain with a $C^{1,\alpha}$ boundary, $\beta,\delta\geq \alpha$ and $p,q$ are like above.
	We show that $\abar$ satisfies \eqref{a}, \eqref{b} and \eqref{c}.  
It holds
	\begin{align}\label{A''pera}
		\sum_{i,j=1}^{n}\partial_ja_1(x,s,\eta)\xi_i\xi_j=&\\
		(\|x\|^\gamma|s|^\delta+1)|\eta|^{p-4}\lg^{q-1}(1+|\eta|)&
		\left[\left((p-2)|\lg(1+|\eta|)+q\frac{|\eta|}{1+|\eta|}\right)\langle\xi,\eta\rangle^2+|\eta|^2\lg(1+|\eta|)|\xi|^2\right]\nonumber 
	\end{align}
	Let now $\mu=\min\{1,p+q-1,p-1\}$. Then
	\begin{align*}
		\left(( p-2)|\lg(1+|\eta|)+q \frac{|\eta|}{1+|\eta|}\right) \langle\xi,\eta\rangle^2+|\eta|^2\lg(1+|\eta|)|\xi|^2\geq\mu |\eta|^{2}|\xi|^2\lg(1+|\eta|)\ \hbox{for}\ \eta,\xi\in \rn\,.
	\end{align*}
	This guarantees that 
	\begin{align*}
		\sum_{i,j=1}^{n}\partial_ja_i(x,s,\eta)\xi_i\xi_j\geq \mu\frac{A'(|\eta|)}{|\eta|}  |\xi|^2\ &\ \hbox{for}\ \eta,\xi\in \rn\,,
	\end{align*}
	and \eqref{a} holds. A simple calculation shows that \eqref{b} holds too.
	For what concerns \eqref{c}, let $M>0$ and take $x,y \in \Omega,\ s,w \in [-M,M],\ \xi \in\rn$. Then 
	\begin{align}\label{a3reg}
		|\abar(x,s,\xi)-\abar(y,w,\xi)|&\leq|\|x\|^\gamma|s|^\delta-\|y\|^\gamma|w|^\delta| A'(|\xi|)\\
		&\leq ((|\|x\|^\gamma-\|y\|^\gamma|)|s|^\delta+\|y\|^\gamma(||s|^\delta-|w|^\delta|)) A'(|\xi|)\nonumber \\
		&\leq	C(\|x-y\|^\alpha+|s-w|^\alpha) A'(|\xi|)\ \hbox{for some}\ C>0\,.\nonumber 
	\end{align}
	For the last inequality in \eqref{a3reg} we used \eqref{in1} and \eqref{in2} further down. They are obtained via the inequalities below.\\
	If $\rho\leq 1$ then there exists $c>0$ such that
	\begin{align}\label{<1}
		|t^\rho-z^\rho|\leq c\frac{|t-z|}{t^{1-\rho}+z^{1-\rho}}\ &\ \hbox{for all}\ t,z>0,
	\end{align}
	If $\rho> 1$ then there exists $c>0$ such that
	\begin{align}\label{>1}
		|t^\rho-z^\rho|\leq c|t-z||(t^{\rho-1}+z^{\rho-1})\ &\ \hbox{for all }\ t,z>0,
	\end{align}
	Take $s,w\in [-M,M]$, for $M>0$. If $0<\delta\leq 1$, from $(|s|^{1-\delta}+|w|^{1-\delta})^{\frac1{1-\delta}}\geq |s|+|w|$ and \eqref{<1}
	\begin{align}\label{in1}
		||s|^\delta-|w|^\delta|&\leq c\frac{||s|-|w||}{|s|^{1-\delta}+|w|^{1-\delta}}\leq c\frac{||s|-|w||^{1-\alpha+\alpha}}{(|s|+|w|)^{1-\delta}} \\
		&\leq c ||s|-|w||^{\alpha}(|s|+|w|)^{\beta -\alpha}\leq c(2M)^{\delta -\alpha}|s-w|^{\alpha}\,.\nonumber
	\end{align}
	A similar argument, when $ 1\leq \delta\leq \alpha$, via \eqref{>1}, leads to
	\begin{align}\label{in2}
		||s|^\delta-|w|^\delta|&\leq c\leq c2^{2 -\alpha}M^{\delta -\alpha}|s-w|^{\alpha}\,.
	\end{align}
	The same holds for $|\|x\|^\gamma-\|y\|^\gamma|,\ x,y\in \Omega$. Thus \eqref{a3} remains true.
\end{example}

For the first existence and regularity Theorem we assume that problem  \eqref{problem} admits a subsolution and a supersolution $\underline{u}$, $\overline{u}\in W^{1,\infty}(\Omega)$.

\begin{theorem}\label{reg}
	Let $\Omega$ be a bounded domain in $\rn$ with a $C^{1,\alpha}$ boundary. 
	Let the functions $A$ and $\abar$ be as in Proposition \ref{Lib}. Assume further that \eqref{a}, \eqref{b} and \eqref{c} hold for all $s\in\r$. Let $\underline{u}$, $\overline{u}\in W^{1,\infty}(\Omega)$ be a subsolution and a supersolution for problem \eqref{problem},  with $\underline{u}(x)<\overline{u}(x)$ a.e. $x\in \Omega$. Let $f:\Omega\times \r\times \rn\to \r$ be a Carath\'{e}odory function satisfying
\begin{align}\label{growthf}
	|f(x,s,\xi)|\leq \sigma (x)+\overline \gamma (s)A'(|\xi|)|\xi|\quad \hbox{for a.e.}\ x\in \Omega,\ \hbox{all}\ s\in [\underline u(x),\overline u(x)],\ \hbox{all}\ \xi\in \rn\,,
	\end{align}
where $\sigma\in L^{\infty}(\Omega)$ and $\overline \gamma:[0,+\infty[\to [0,+\infty[$ is locally essentially bounded.\\	
Then problem \eqref{problem} admits at least a solution  $u\in C_0^{1,\beta} (\overline \Omega)$. Moreover $\underline u(x)\leq u(x)\leq \overline u(x)$ a.e in $\Omega$.
\end{theorem}
{\bf Proof}. Let $M=\max\{\| \overline u\|_\infty,\| \underline u\|_\infty\}+1>0$,  $R>\max\{\|\nabla \overline u\|_\infty,\|\nabla \underline u\|_\infty\}$. Let $I=[-M,M]$ and $\overline \gamma=\|\overline \gamma\|_{L^\infty(I)}$. Then $\overline \gamma<+\infty$ and $\overline \gamma(s)\leq \overline \gamma$ for a.e. $s\in I$. Let $I_0\subset I$ be a set of null measure, such that $\overline \gamma(s)>\overline \gamma$ for all $s\in I_0$. For $s\in I_0$ it holds
\begin{align*}
	|f(x,s,\xi)|=\lim_{t\to s}|f(x,t,\xi)|=\liminf_{t\to s}|f(x,t,\xi)|& \leq \sigma(x)+\liminf_{t\to s}\overline \gamma(t)A'(|\xi|)|\xi|\\
	\leq \sigma(x)+\overline \gamma A'(|\xi|)|\xi| &\quad \hbox{for a.e.} \; x\in \Omega, \; \hbox{all} \ \xi\in\rn\,.\nonumber
\end{align*}
Let the truncated functions $a_i^M$ and $f_R$ be defined as
\begin{equation*}\label{aM}
	a_i^M(x,s,\xi)=\left\{\begin{array}{cc}
		a_i(x,-M,\xi) & \hbox{if}\ s\leq -M,\\
		a_i(x,s,\xi) & \hbox{if}\ -M<s<M \,,\\
			a_i(x,M,\xi) & \hbox{if}\ s\geq M\,,
	\end{array}
	\right.
\end{equation*}
and
\begin{equation*}\label{fR}
	f_R(x,s,\xi)=\left\{\begin{array}{cc}
		f(x,s,\xi) & \hbox{if}\ |\xi|\leq R,\\
		f(x,s,\xi)\cdot\frac{A'(R)R}{A'(|\xi|)|\xi|} & \hbox{if}\ |\xi|> R \,.
	\end{array}
	\right.
\end{equation*}
Consider the problem
\begin{equation}\label{problem1}
	\begin{cases} 
		-div(\abar^M(x,u,\nabla u))=f_R(x,u,\nabla u) & {\rm in}\,\,\,\Omega \\
		u  =0 & {\rm on}\,\,\,\partial \Omega \,.
	\end{cases}
\end{equation}
In view of the choice of $R$, $\underline u$ and $\overline u$ are a subsolution and a supersolution to \eqref{problem1} respectively. The function $f$ satisfies satisfies the hypotheses of Proposition \ref{Lib}. In particular, \eqref{d} holds ith $\Lambda_1=\max\{\|\sigma\|_\infty, \overline \gamma\}$. 
Since $|f_R|\leq |f|$ the same holds for $f_R$ whatever $R$ is. Also $\abar^M$ satisfies \eqref{a}, \eqref{b} and \eqref{c}.
Due to Proposition \ref{Lib} there exist two positive constants $0<\beta\leq 1$ and $C$,  independent from $R$, such that any solution to \eqref{problem1} belongs in $C_0^{1,\beta}(\overline \Omega)$ and $\|u\|_{C_0^{1,\beta}(\overline \Omega)}\leq C$. Choosing $R>C$ we deduce that $u$ is a solution to \eqref{problem}.
\qed

The next Theorem is related to Theorem \ref{Teo1}, because it deals with problem \eqref{problemvar}. Now, we are concerned with the existence of a regular solution to  \eqref{problemvar}. 
\begin{theorem}\label{main1r}
Let $\Omega$ be a bounded domain in $\rn$ with $C^{1,\alpha}$ boundary. 
Let $A$ be a Young function satisfying \eqref{dg}, and let $\abar:\Omega\times\rn\to \rn$ be a vector valued function satisfying \eqref{a} and \eqref{b} for a.e. $x\in\Omega$, all $\xi\in\rn$, and such that $a_i(x,0)\equiv 0$ for a.e. $x\in\Omega$, all $i=1,\ldots,n$. We further assume that \eqref{a4} holds. Let $f:\Omega\times \mathbb{R}\times \mathbb{R}^n\to \mathbb{R}$ be a Carath\'{e}odory function fulfilling
	\begin{align}\label{growth f2+}
		-\rho_2(x)-g_2(s)-\overline \gamma(s)A'(|\xi|)|\xi|&\leq f(x,s,\xi)\leq \rho_1(x)+g_1(s)\ \hbox{for a.e.}\ x\in \Omega,\ \hbox{all}\ s\geq 0,\\ 
		\ \hbox{all}\ \xi\in\rn,\  f(x,0,0)\geq 0 &\ \hbox{in}\ \Omega,\ \hbox{and}\  f(x,0,0)>0\ \hbox{on a set of positive measure},\nonumber
	\end{align}
or 
\begin{align}\label{growth f2-}
-\rho_1(x)-g_1(|s|)	\leq f(x,s,\xi)&\leq \rho_2(x)+g_2(|s|)+\overline \gamma(s)A'(|\xi|)|\xi|\	\hbox{for a.e.}\, x\in \Omega,\ \hbox{all}\ s\leq 0,\\ 
\ \hbox{all}\ \xi\in\rn,\ f(x,0,0)\leq 0 &\ \hbox{in}\ \Omega,\hbox{and}\ f(x,0,0)<0\ \hbox{on a set of positive measure},\nonumber
\end{align}
Here $\rho_1,\rho_2:\Omega\to [0,+\infty[$ are two measurable functions, $\rho_1,\,\rho_2\in L^\infty(\Omega)$; $g_1$ is like in Theorem \ref{Teo1}, $g_2:[0,+\infty[\to [0,+\infty[$ is a  non-decreasing function such that $g_2(0)=0$ and $\overline \gamma (s)$ is a locally essentially bounded function.\\
Then problem $(P)$ has a nontrivial, solution $u\in C_0^{1,\beta}(\overline \Omega)$.\\
If \eqref{growth f2+} jolds, then $u\geq 0$ in $\Omega$. In the other case $u\leq 0$ in $\Omega$.
\end{theorem}




{\bf Proof.}
From the proof of Theorem \ref{Teo1} we know that, when \eqref{growth f2+} is in force, there exists a nontrivial solution $\overline{u}\geq 0$, to problem \eqref{problemvar}.\\
Remark \ref{rembound} guarantees that $\overline{u}$ is bounded. From Proposition \ref{Lib}, we have that $\overline{u}\in C_0^{1,\beta}(\overline{\Omega})$.
The inequalities in \eqref{growth f2+} show that $\overline{u}$ is a supersolution to problem \eqref{problem} and $\underline{u}=0$ is a subsolution to problem \eqref{problem}. The assumptions on $\rho_2$ guarantee that $\underline{u}=0$ is not a solution. If we put $M=\|\overline u\|_\infty$, $\overline \gamma=\|\overline\gamma\|_{L^\infty(0,M)}$, $\sigma(x)=\max\{\rho_1(x)+g_1(M) ,\; \rho_2(x)+g_2(M)\}$ for a.e. $x\in \Omega$, then \eqref{growth f2+} leads to
\begin{align}
|f(x,s,\xi)|\leq \sigma(x)+\overline\gamma A'(|\xi|)|\xi|&\quad \hbox{for a.e.} \; x\in \Omega, \; s\in [0,\overline u(x)],\; \xi\in \mathbb{R}^n\,.
\end{align}
So, from Theorem \ref{reg}, problem \eqref{problemvar} admits at least a nontrivial solution $u\in C_0^{1,\beta}(\overline{\Omega})$ such that $0\leq u\leq \overline{u}$.\\
For the other case, it is enough to put $f_1(x,s,\xi)=-f(x,-s,-\xi)$ and to use the first part of the proof.\qed

\section{Examples}\label{sec5}

In this Section we present some applications of Theorem \ref{Tnew} (Subsection \ref{ex1}) where the structure of $f$ and the condition $\abar(x,s,0)\equiv 0_{\rn}$ guarantee the existence of constant sub and supersolutions for \eqref{problem}. This situation is really interesting and meaningful because hilights how in this setting the growth of $\abar$ with respect to $s$ can be whatever we want. Then we present an application of Theorem \ref{Teo1} (Subsection \ref{ex2}) and finally some applications of Theorem \ref{main1r} (Subsection \ref{ex3}). In the latter case the growth of $\abar$ with respect to $s$ comes once again into play (see equation \eqref{c}).\\
In all the examples $\Omega$ is an open subset in $\rn$ with finite measure.


\subsection{Applications of Theorem \ref{Tnew}}\label{ex1}

In this example we do not impose growth conditions with respect to $s$ for $\abar$, but only with respect to $x$ and $\xi$. 
The structure of the convection term guarantees that \eqref{problem} has a pair of constant sub and supersolutions.\\
Let $0<\delta<p-1$, $p+q>1$, $a:\Omega\to [0,+\infty[$ be a measurable function, $a\in L^{\frac{p}{\delta}}(\Omega)$, and let $b:[0,+\infty[\to [0,+\infty[$ be a
 continuous function.\\
 Consider the problem
\begin{equation}\label{xsuf}
	\begin{cases}
	- {\rm div}\left(\left(a(x)b(|u|)|\nabla u|^{p-2-\delta}\lg^{q(1-\frac\delta{p-1})}(1+|\nabla u|) +|\nabla u|^{p-2}\lg^q(1+|\nabla u|)\right)\nabla u\right) =f(x,u,\nabla u) & {\rm in}\,\, \Omega \\
		u  =0 & {\rm on}\,\partial \Omega \,,
	\end{cases}
\end{equation}
The Young function $A$ governing the differential operator $\abar$ obeys \eqref{Aex}.\\
The convection term $f:\Omega\times\r\times\rn\to \r$ is defined as
\begin{align}\label{fpq}
	f(x,s,\xi)=(h(x)+k(|\xi|))g(s)\quad  \hbox{for}\ (x,s,\xi)\in \Omega \times \r\times \rn\,. \end{align} 
Here $h:\Omega\to [0,+\infty[$ is a measurable function, $h\in L^{\widetilde A_n}(\Omega)$, 
$k:[0,+\infty]\to \r$ is a continuous function, $k(0)>0$, and $k$ has the following behavior near infinity 
\begin{align}\label{h1pq}
	\left\{\begin{array}{ll}
	|k(s)|\approx s^{\frac{p}{(p^*)'}}\lg^{r}(s) &\ \hbox{for some}\ r<\frac{(n+1)q}{n},\ \hbox{when}\ p<n\,,\\
	|k(s)|\approx s^{n}\lg^{r}(s)&\ \hbox{for some},\ r<q-1+\frac{q+1}{n},\ \hbox{when}\ p=n,\ q<n-1\,,\\
	|k(s)|\approx s^{n}\lg^{n-1}(s)\lg^{-\frac 1r}(\lg s)&\ \hbox{for some},\ r<\frac{n-1}{n},\ \hbox{when}\ p=n,\ q=n-1\,,\\ 
|k(s)|\approx s^{p}\lg^{r}(s)&\ \hbox{for some},\ r<q,\ \hbox{when}\ p>n\ \hbox{or}\ p=n, q>n-1\,.\end{array} \right.
\end{align}
The function $g:\r\to\r $ is continuous, $g(s)>0$ for $s\in [0,\overline s)$ and $g(\overline s)=0$.
First of all we note that $u_1=0$ and $u_2=\overline s$ are a subsolution and a supersolution to \eqref{xsuf} and $u\equiv 0$ is not a solution. The growth of the functions $h$ and $k$ quarantee that $f$ satisfies \eqref{growth f'}.\\
Due to the continuity of $b$, and taking into account 
Example \ref{exA} (with $\beta=\beta_1=0$), we see that conditions \eqref{a1}, \eqref{a2} and \eqref{a3bis} hold, for a.e. $x\in\Omega$, all $s\in [0,\overline s]$,\ all $\xi\in \rn$.
Thus, by Theorem \ref{Tnew}, problem \eqref{xsuf} has a nontrivial solution $u\in [u_1,u_2]$.\\
The same arguments work for different choices of $h$ and $g$. In particular we see that Theorem \ref{Tnew} works well with all nonlinearities having two zeros $s_1$ and $s_2$, with $s_1<0<s_2$, ($f(x,s_1,0)=f(x,s_2,0)=0$ for all $x\in\Omega$) and $f(x,s,0)$ has  constant sign for $s\in ]s_1,s_2[$, or for which  $f(x,0,0)$ has constant sign and there exists $s\in\r$ such that $f(x,s,0)\equiv 0$ and $s\cdot f(x,0,0)>0$. 

\subsection{Applications of Theorem \ref{Teo1}\label{ex2}}
Let $p,q,r\in \r$ be such that $1<r<p<n$, $1<q<p$, and let $m<0$.\\ 
Let $a,\rho:\Omega \to [0,+\infty[$ be two measurable functions, $a\in L^{\frac p{p-r}}(\Omega)$ and $\rho>0$ on a subset of $\Omega$ having positive measure.\\ We show that problem 
\begin{equation}\label{xuf}
	\begin{cases}
		- {\rm div}\left(\left(a(x)|\nabla u|^{r-2} +|\nabla u|^{p-2}\right)\nabla u\right) =\frac{\rho(x)+|u|^{q-1}}{1+|\nabla u|}-|\nabla u|^{\frac{p}{{p^*}'}}|\lg(|\nabla u|)|^{\frac{m}p^*} & {\rm in}\,\, \Omega \\
		u  =0 & {\rm on}\,\partial \Omega \,,
	\end{cases}
\end{equation}
has a nontrivial solution $u\in W_0^{1,p}(\Omega)$, $u\geq 0$.\\
The functions $\Phi(x,\xi)=a(x)\frac{|\xi|^{r}}{r} +\frac{|\xi|^{p}}{p}$, and  $\abar(x,\xi)=\left(a(x)|\xi|^{r-2} +|\xi|^{p-2}\right)\xi$ satisfy \eqref{a_1'}, \eqref{a_2'}, \eqref{a_3'} and \eqref{a4}. Condition \eqref{growth f2} holds too, with $\rho_1(x)\equiv \rho(x)$, $g_1(|s|)=|s|^{q-1}$, $\rho_2(x)=-M^{\frac{p}{{p^*}'}}|\lg(|M|)|^{\frac{m}p^*}$, for a suitable constant $M>0$, $g_2(|s|)\equiv 0$ and $E(t)\approx t^{p^*}\lg^m(t)$, for $t>>1$.\\Thus from Theorem \ref{Teo1} problem \eqref{xuf} has a nontrivial solution $u\in W_0^{1,p}(\Omega)$, $u\geq 0$.

\subsection{Applications of Theorem \ref{reg}\label{ex3}}
Let $\Omega$ be a bounded domain with a $C^{1,\alpha}$ boundary. Let $\gamma\geq \alpha$, $p>1$ and $q\in \r$, satisfying $p+q-1>0$. Let $h:\Omega \to [0,+\infty[$ be a measurable function, $h\in L^{\infty }(\Omega)$, $h\neq 0$ on a subset of $\Omega$ having finite measure, and let $g:\r\to \r$ be a continuous function such that $g(s)>0$ for all $s\in [0,\overline s[$ and $g(\overline s)\equiv 0$. 
We consider the problem 
\begin{equation}\label{xuxifreg}
	\begin{cases}
		- {\rm div}\left(\|x\|^\gamma e^{|u|}|\nabla u|^{p-2} \lg^q(1+|u|)\nabla u\right) =\left(h(x)+|\nabla u|^p\lg^{q}(1+|\nabla u|)\right)g(u)  & {\rm in}\,\, \Omega \\
		u  =0 & {\rm on}\,\partial \Omega \,,
	\end{cases}
\end{equation}
and prove that it has a nontrivial solution $u\in W_0^{1,A}(\Omega)$, $u\geq 0$. The Young function $A$ is defined via \eqref{A'ex}.\\
Using the Mac Laurin expansion of $k(t)=e^t$ we see that $|e^{|s|}-e^{|w|}|\leq \frac{e^M(e^M-1)}{M^\alpha}|s-w|^\alpha$ for all $s,w\in [-M,M]$. Thus, following the arguments used in the Example \ref{abarforreg}, we can prove that the operator $\abar$ and the function $f$, defined respectively as $\abar(x,s,\xi)=\|x\|^\gamma e^{|s|}|\xi|^{p-2} \lg^q(1+|\xi|)\xi$  and $f(x,s,\xi)=\left(h(x)+|\xi|^p\lg^{q}(1+|\xi|)\right)g(s)$ satisfy the hypotheses of Theorem \ref{reg}. Furthermore, $\underline u\equiv 0$ and $\overline u=\overline s$ are a subsolution and a supersolution to \eqref{xuxifreg}, and $\underline u\equiv 0$ is not a solution. Thus, from Theorem \ref{reg}, the problem which we are dealing with, has a regular solution $u\in [0,\overline s]$.

\end{document}